\pdfoutput=1
\title{A Plane Wave Discontinuous Galerkin Method with a
Dirichlet-to-Neumann Boundary Condition for the
Scattering Problem in Acoustics}
\author{
        Shelvean Kapita \\
                Institute for Mathematics and Its Applications\\
        University of Minnesota, 207 Church St SE\\
        Minneapolis, Minnesota 55455
            \and
        Peter Monk\\
        Department of Mathematical Sciences\\
        University of Delaware\\
        Newark, Delaware 19716
}
\date{\today}

\documentclass[10pt]{article}

\usepackage{color} 
\usepackage{graphicx}
\usepackage{bbm}
\usepackage{amsmath}
\usepackage{amssymb}
\usepackage{amsfonts}   
\usepackage{mathtools}
\usepackage{multirow}
\usepackage{array}
\newcommand{\jmp}[1]{\left[\!\left[#1\right]\!\right]}                     
\newcommand{\avg}[1]{\left\{\!\!\left\{#1\right\}\!\!\right\}}                   
\usepackage{pgfplots}
\usepackage{float}
\usepackage[bookmarks=false]{hyperref}
\newcommand{\cE}{\mathscr{E}}
\usepackage{hyperref}
\newcommand{\cEI}{\mathscr{E}_{\!_I}}
\newcommand{\cED}{{\mathscr{E}}_{\!_D}}
\newcommand{\cEB}{{\mathscr{E}}_{\!_R}}

\newcommand{\OmD}{\Omega_{\!_D}}

\newtheorem{theorem}{Theorem}[section]
\newtheorem{lemma}[theorem]{Lemma}
\newtheorem{proposition}[theorem]{Proposition}
\newtheorem{corollary}[theorem]{Corollary}

\newcommand{\diff}[2]{\frac{\partial #1}{\partial #2}}

\newcommand{\Imag}{\operatorname{Im}}
\newcommand{\ThN}{\mathcal{S}_{\!_N}}

\newcommand{\hf}{\frac{1}{2}}
\newcommand{\htf}{\frac{3}{2}}

\newcommand{\bx}{\textbf{\textit{x}}}
\newcommand{\bdr}{\textbf{\textit{d}}}
\newcommand{\bn}{\textbf{\textit{n}}}
\newcommand{\DtN}{\mathcal{S}}
\newcommand{\Pro}{\mathcal{P}}
\newcommand{\bsig}{\boldsymbol{\sigma}}

\newcommand{\pwdg}[1]{{PW}(#1)}
\newcommand{\ep}{\scalebox{1.0}{$\square$}}

\definecolor{lightgray}{gray}{0.75}
\newcommand{\un}{u^{\!_N}}
\usepackage{mathrsfs} 

\begin{document}

\maketitle

\tableofcontents

\begin{abstract}
We consider the numerical solution of an acoustic scattering problem by the Plane Wave Discontinuous Galerkin Method (PWDG) in the exterior of a bounded domain in $\mathbb{R}^2$. In order to apply the PWDG method, we introduce an artificial boundary to truncate the domain, and we impose a non-local Dirichlet-to-Neumann (DtN) boundary conditions on the artificial curve. To define the method, we introduce new consistent numerical fluxes that incorporate the truncated series of the DtN map. Error estimates with respect to the truncation order of the DtN map, and with respect to mesh width are derived. Numerical results suggest that the accuracy of the PWDG method for the acoustic scattering problem can be significantly improved by using DtN boundary conditions. 
\end{abstract}


%

%
\section{Introduction}
\label{introduction.sec}
Acoustic, elastic and electromagnetic scattering problems arise in many areas of physical and engineering interest, in areas as diverse as radar, sonar, building acoustics, medical and seismic imaging. Mathematically, the problem of acoustic scattering is often modeled by the Helmholtz equation in the unbounded region exterior to a bounded obstacle $\Omega_{\!_D}\subset\mathbb{R}^2$. The main difficulties in the numerical simulation of the scattering problem arise from the unbounded nature of the domain, and from the usually oscillatory nature of the solutions. Recent work has led to the development of algorithms that are better able to handle the highly oscillatory nature of the solutions at high frequency. These algorithms incorporate in their trial and test spaces oscillatory functions such as plane waves, Fourier Bessel functions, and products of low order polynomials with plane waves in order to capture qualitative information about the oscillatory character of the solution. These methods include the Partition of Unity Method (PUM) of Melenk and Babuska,\cite{babuskamelenk} the Discontinuous Enrichment Method (DEM) \cite{dem}, the Ultra Weak Variational Formulation (UWVF) of Cessenat and Despr{\' e}s \cite{cd}. The UWVF has been applied to the Maxwell equations \cite{hmm}, linear elasticity \cite{lhm}, acoustic fluid-solid interaction \cite{hkm} and to thin clamped plate problems \cite{lhm1}. More recently, the Plane Wave Discontinuous Galerkin (PWDG) method has been studied by Hiptmair et al. \cite{HMP11, HMP13, HMP15, GH14} as a generalization of the UWVF method. In this case the problem is cast in the form of a Discontinuous Galerkin method, with possibly mesh dependent penalty parameters. 

Hitherto, since the focus has been to understand plane wave discretization, only approximate first order absorbing boundary conditions have been considered in PWDG methods to truncate the exterior domain. However, unless the computational domain is large, first order absorbing boundary conditions result in large errors due to reflections from the artificial boundary. To reduce errors due to boundary reflections, without taking the absorbing boundary far from the scatterer, several numerical techniques have been considered. These include high order absorbing boundary conditions \cite{Higdon}, the Dirichlet-to-Neumann (DtN) mapping \cite{Hsiao, kmdtn, kdtn}, and perfectly matched layers \cite{berenger}. The DtN map can be enforced via boundary integral equations or Fourier series expansions resulting from the method of separation of variables. In this paper, we will couple the Fourier series DtN method with PWDG.

In \cite{kmdtn, kdtn}, Koyama derives error estimates for the DtN finite element method, considering both the error due to truncation of the DtN series and due to finite element discretization. Hsiao et. al \cite{Hsiao} derive error estimates for the conforming finite element DtN method, and present numerical results that show optimal convergence in the $L^2$ and $H^1$ norms using conforming piecewise linear finite elements. Following \cite{kmdtn, kdtn, Hsiao}, we consider an approximate boundary value problem with a truncated DtN series. Our main contribution in this work is to choose new numerical fluxes on the artificial boundary that enforce the truncated DtN boundary condition, and that are consistent with respect to the solution of the truncated boundary value problem. In the Ph.D. thesis \cite{moiolathesis}, A. Moiola derives approximation results for solutions of the Helmholtz equation by plane waves, using an impedance first order absorbing boundary condition. We modify these arguments to take into account the non local DtN series, and to derive error estimates in this case. Since the DtN map naturally maps from $H^\hf$ on the artificial boundary, but plane wave basis functions are only in $L^2$, our error estimates deteriorate with order $N$ as $N\rightarrow\infty$. This problem can be overcome by using the Neumann-to-Dirichlet map (NtD), and is part of ongoing work to incorporate generalized impedance boundary conditions into the PWDG method.

This paper is organized as follows: In Section~\ref{notation.sec}, we consider the continuous problems with the full and truncated DtN series, and prove stability estimates in both cases. In Section~\ref{pwdg.sec}, we introduce the PWDG method for solving the scattering problem, and show the well-posedness and consistency of the PWDG scheme. In Section~\ref{error.sec}, we modify the error estimates in \cite{moiolathesis} to take into account the DtN boundary condition. Section~\ref{implm.sec}, we describe the numerical implementation of the DtN map via a projection into the space of trigonometric polynomials. In Section~\ref{results.sec}, numerical results are presented to show the effect of the wavenumber, truncation of the DtN series, number of plane waves, and mesh size on the convergence of the numerical scheme. We end in Section~\ref{conclusion.sec} with some concluding remarks.

\section{The Scattering Problem}
\label{notation.sec}
We consider the scattering of an acoustic wave from a sound-soft impenetrable obstacle occupying a region $\OmD\subset\mathbb{R}^2$. We assume that the boundary of $\Omega_{\!_D}$, denoted $\Gamma_{\!_D}$ is smooth, and the  exterior $\mathbb{R}^2\backslash\Omega_{\!_D}$ is connected. The scattered field $u$ satisfies the following Dirichlet boundary value problem in $\mathbb{R}^2\backslash\overline{\OmD}$:
\begin{subequations}\label{helmholtzbc}
\begin{align}
&\Delta u + k^2 u= 0\;\;\mbox{in}\;\mathbb{R}^2\backslash\overline{\Omega}_{\!_D} \label{helmeqn}\\
& u = g\;\;\mbox{on}\;\Gamma_{\!_D}:=\partial\Omega_{\!_D}\label{dirbc}\\
& \lim_{r\rightarrow\infty}r^\hf\left(\diff{u}{r}+iku \right) = 0,\label{src}
\end{align}
\end{subequations}
where $r = |\bx|$, $\bx = (x_1,x_2)$, and $g=-u^{\tiny{inc}}$, where $u^{\tiny{inc}}$ is a known incident field. We will assume that the wavenumber $k>0$ is constant, although the algorithm applies to the more general cases where $k$ is piecewise constant in each element. The well-posedness of the exterior Dirichlet problem~(\ref{helmholtzbc}) is demonstrated for example in Theorem  3.11 of \cite{coltonkress}.

To apply a domain based discretization, it is necessary to introduce an artificial domain $\Omega_{\!_R}$ with boundary $\Gamma_{\!_R}$ enclosing $\Omega_{\!_D}$ such that $\mbox{dist}(\Gamma_{\!_R},\Gamma_{\!_D})>0,$ and to introduce suitable boundary conditions on $\Gamma_{\!_R}$ that take into account wave propagation in the infinite exterior of $\Omega_{\!_R}$. The boundary value problem (\ref{helmholtzbc}) is then replaced by a boundary value problem posed in the annulus bounded by $\Gamma_{\!_R}$ on the outside and $\Gamma_{\!_D}$ on the inside.

To date, a crucial step taken in PWDG methods for the Helmholtz equation is to use the \emph{first order absorbing boundary condition} 
\begin{equation}
\begin{aligned}\label{impedancebc}
&\frac{\partial u}{\partial \bn}+iku = 0\;\mbox{on}\;\Gamma_{\!_R}
\end{aligned}
\end{equation}
In this paper, as in \cite{wildemonk, Hsiao}, we are going to use the Dirichlet-to-Neumann map now in the context of the PWDG method.

\subsection{Non-local Boundary Value Problem}
The scattering problem (\ref{helmholtzbc}) is equivalent to the following boundary value problem (see, e.g. Johnson and Nedelec \cite{johnsonnedelec}):
\begin{equation}\label{unbounded}
\left. \begin{aligned}
&\Delta u +k^2 u = 0\;\;\mbox{in}\;\Omega,\;
u = g\;\;\mbox{on}\;\Gamma_{\!_D},\\
&u = w,\;\;\frac{\partial u}{\partial\bn} = \frac{\partial w}{\partial\bn}\;\;\mbox{on}\;\Gamma_{\!_R},\\
&\Delta w + k^2 w = 0\;\;\mbox{in}\;\mathbb{R}^2\backslash\overline{\Omega},\\
&\lim_{r\rightarrow\infty}r^\hf\left(\frac{\partial w}{\partial r} +ikw\right)=0.
\end{aligned}
\right\}
\end{equation}
Here $u,w$ are the scattered fields in the interior and exterior of $\Omega$ respectively, and $g\in H^\hf(\Gamma_{\!_R})$. If $w$ is known on $\Gamma_{\!_R}$, the normal derivative ${\partial_\bn w}$ can be computed by solving for $w$ in $\mathbb{R}^2\backslash\overline{D}$. The DtN map $\DtN:H^\hf(\Gamma_{\!_R})\rightarrow H^{-\hf}(\Gamma_{\!_R})$ is defined as 
$$\DtN:w|_{\Gamma_{\!_R}}\rightarrow \partial_\bn w|_{\Gamma_{\!_R}}.$$

If the truncating boundary $\Gamma_{\!_R}$ is a circle, the map $\DtN$ can be written explicitly as a series involving Hankel functions. By using the polar coordinate system, separation of variables shows that the general solution of the homogeneous Helmholtz equation $\Delta w + k^2 w=0\;\;\mbox{in $\mathbb{R}^2\backslash\overline{\Omega}$}$ is 
\begin{eqnarray}\label{seriessoln}
w(r,\theta)&=&\sum_{m\in\mathbb{Z}}\left[\alpha_{\!_m}H_m^{(1)}(kr)+\beta_{\!_m}H_m^{(2)}(kr)\right]e^{im\theta},
\end{eqnarray}
where $H^{(1,2)}_m(z)$ are Hankel functions of first and second kind, and of order $m\in\mathbb{Z}$. The Hankel functions are in turn defined by Bessel functions $J_m(z)$ and Neumann functions $Y_m(z)$
$$H^{(1,2)}_m(z) = J_m(z)\pm iY_m(z).$$
For an introduction to Bessel and Hankel functions in the context of the Helmholtz equation, see Colton and Kress \cite{coltonkress}, or Cakoni and Colton \cite{cakonicolton}. 

Only the Hankel functions of the \emph{second kind} are consistent with the Sommerfeld  radiation condition (\ref{src}), so only solutions of the form $H^{(2)}_m(kr)e^{im\theta}$ represent outgoing waves. This implies that the coefficients $\alpha_{m}$ in the series expansion (\ref{seriessoln}) vanish.
If $w(R,\theta)\in H^\hf(\Gamma_{\!_R})$ is given, then we can write $w$ as a Fourier series,
\begin{eqnarray*}
w(R,\theta) &=& \sum_{m\in\mathbb{Z}}w_{\!_m}e^{im\theta}
\end{eqnarray*}
where the Fourier coefficients $w_{\!_m}$ are given by 
\begin{eqnarray*}
w_{\!_m}(R)&=& \frac{1}{2\pi R}\int_{\Gamma_{\!_R}}w(R,\varphi)e^{-im\varphi}\;d\varphi.
\end{eqnarray*} 
Thus, the solution $w$ of the Helmholtz problem for $r\geq R$ is
\begin{eqnarray}\label{seriessoln1}
w(r,\theta) &=& \sum_{m\in\mathbb{Z}}w_{\!_m}(R)\frac{H^{(2)}_m(kr)}{H^{(2)}_m(kR)}e^{im\theta}.
\end{eqnarray}
Taking the normal derivative of $w(r,\theta)$, which is simply the radial derivative $\partial w/\partial r$ we can write an explicit form of the DtN map $\DtN$:
\begin{eqnarray}\label{dtndfn}
\DtN w(R,\theta):=\frac{\partial w}{\partial r}(R,\theta) &=& \sum_{m\in\mathbb{Z}}k\frac{H^{(2)'}_m(kr)}{H^{(2)}_m(kR)}w_m(R)\;e^{im\theta}.
\end{eqnarray}
Using the DtN map, we may restrict the domain of problem (\ref{unbounded}) to $\Omega$ and the equations for $u$ become
\begin{equation}\label{dtnbvp}
\left. \begin{aligned}
&\Delta u + k^2 u = 0,\;\;\mbox{in}\;\Omega\\
& u = g,\;\;\mbox{on}\;\Gamma_{\!_D}\\
& \frac{\partial u}{\partial\bn}-\DtN u=0\;\;\mbox{on}\;\Gamma_{\!_R}.
\end{aligned}
\right\}
\qquad
\end{equation}
The boundary condition $\partial u/\partial\bn-\DtN u=0$ on $\Gamma_{\!_R}$ is an \emph{exact} representation of wave propagation in the exterior domain. The solution $u$ of the boundary value problem~(\ref{dtnbvp}) is the restriction to $\Omega$ of the unique solution of the exterior scattering problem~(\ref{helmholtzbc}). 

\subsection{Stability of the adjoint problem}
 Since we will later use duality arguments to derive error estimates, we consider the following adjoint problem to~(\ref{dtnbvp})
 \begin{equation}\label{adjointdtnbvp}
\left. \begin{aligned}
&\Delta u + k^2 u = f,\;\;\mbox{in}\;\Omega\\
& u = 0,\;\;\mbox{on}\;\Gamma_{\!_D}\\
& \frac{\partial u}{\partial\bn}-\DtN^{\star} u=0\;\;\mbox{on}\;\Gamma_{\!_R}.
\end{aligned}
\right\}
\qquad
\end{equation} 
where $\DtN^{\star}$ is the $L^2$ adjoint of the DtN map on $\Gamma_{\!_R}$ and $f\in L^2(\Omega)$ is such that $\mbox{supp}(f)\subset\Omega$. Wavenumber dependent stability bounds for the scattering problem with the impedance boundary condition $\partial u/\partial\bn-iku=g_{\!_R}$ on $\Gamma_{\!_R}$ are shown in \cite{HMP13}, Theorem 2.2.
In the following theorem, we state stability bounds for the scattering problem with a DtN map on $\Gamma_{\!_R}$. We will assume that the boundary of the scatterer is smooth, so that $u\in H^2(\Omega)$. Although not needed here, the following result is of independent interest.
\begin{proposition}\label{contdepn4}
Let $u\in H^2(\Omega)$ be the analytical solution of the adjoint problem~(\ref{adjointdtnbvp}). Then there exist positive numbers $C_1$ and $C_2$ independent of $u$ and $f$, but whose dependence on $k$ and $R$ is known such that
\begin{eqnarray}
&&|u|_{1,\Omega}+k\|u\|_{0,\Omega}\leq C_1\|f\|_{\!_{0,\Omega}},\label{stab1}\\
&&|u|_{2,\Omega}\leq  C_2\|f\|_{\!_{0,\Omega}}\label{stab2}.
\end{eqnarray}
\end{proposition}
{\bf Proof}:
The solution $u$ of the adjoint problem~(\ref{adjointdtnbvp}) can be extended analytically by Hankel functions of the first kind to the exterior region $\mathbb{R}^2\backslash\overline{\Omega}$. Denote still by $u$ this analytic extension in the region $\mathbb{R}^2\backslash\overline{D}$. Let $\widetilde{\Omega}:=B_{\!_{2R}}({\bf 0})\backslash\overline{D}$ be the annulus bounded by the circle $\Gamma_{\!_{2R}}$ on the outside and by $\Gamma_{\!_D}$ on the inside, where $\Gamma_{\!_{2R}}$ is a circle of radius $2R$ centered at the origin.

Let $\widetilde{u} = \chi u$ where $\chi\in\mathcal{C}_0^\infty(\widetilde{\Omega})$,  $0\leq\chi\leq 1$, is a smooth cut-off function equal to one in a neighborhood of $\overline{\Omega}$ and zero in a neighborhood of $\Gamma_{\!_{2R}}$. Then $\widetilde{u}$ satisfies 
\begin{eqnarray*}
&&\Delta\widetilde{u}+k^2\widetilde{u} = \widetilde{f}\;\;\mbox{in}\;\;\widetilde{\Omega},\\
&& \widetilde{u}=0\;\;\mbox{on}\;\;\Gamma_{\!_D},\\
&& \partial_\bn\widetilde{u}+ik\widetilde{u}=0\;\;\mbox{on}\;\;\Gamma_{\!_{2R}},
\end{eqnarray*}
where \begin{equation*}
\widetilde{f}=\left\lbrace\begin{aligned} & f\;\;\mbox{in}\;\;\Omega\\
& \Delta(\chi u)+ k^2\chi u.\;\;\mbox{in}\;\;\widetilde{\Omega}\backslash\overline{\Omega}
\end{aligned}
\right.
\end{equation*}

By the first stability estimate in inequalities (3.5) of \cite{HMP11}, there exists some constant $C$ independent of $k,\widetilde{u}$ and $\widetilde{f}$ such that
\begin{eqnarray}\label{stabres1}
|\widetilde{u}|_{{1,\widetilde{\Omega}}}+k\|\widetilde{u}\|_{{0,\widetilde{\Omega}}}\leq C R\|\widetilde{f}\|_{\!_{0,\widetilde{\Omega}}}.
\end{eqnarray} 
The product rule shows that $\Delta(\chi u) = \chi\Delta u + u\Delta\chi +2\nabla\chi\cdot\nabla u$. Hence $$\widetilde{f} = \chi f +  u\Delta\chi +2\nabla\chi\cdot\nabla u$$ where $f$ is extended by zero to the exterior of $\Omega$. Since we can choose $$|\chi|\leq 1,\;\;\; |\nabla\chi|\leq C/R, \;\;\;|\Delta\chi|\leq C/R^2 $$ we have that 
$$\|\widetilde{f}\|_{{0,\widetilde{\Omega}}}\leq C_0\left(1+\frac{1}{kR}+\frac{1}{k^2R^2}\right)k\left(|u|_{{1,\widetilde{\Omega}}}+k\|u\|_{{0,\widetilde{\Omega}}}\right).$$
By Lemma~3.5 of~\cite{wildemonk}, the solution ${u}$ of~(\ref{adjointdtnbvp}) satisfies the stability bound $$k\left(|{u}|_{1,\widetilde{\Omega}}+ k\|{u}\|_{0,\widetilde{\Omega}}\right)\leq \left(1+4\sqrt{2}kR\right)\|f\|_{{0,\widetilde{\Omega}}}=\left(1+4\sqrt{2}kR\right)\|f\|_{{0,{\Omega}}}$$ Then we have 
\begin{eqnarray*}
|u|_{1,\Omega}+k\|u\|_{0,\Omega}&\leq& |\widetilde{u}|_{1,\widetilde{\Omega}}+k|\widetilde{u}|_{0,\widetilde{\Omega}}\\
&\leq& C_0R\|\widetilde{f}\|_{0,\widetilde{\Omega}}  \\
&\leq& C_1\|f\|_{{0,\Omega}}
\end{eqnarray*}
where $$C_1:= C_0R\left(1+4\sqrt{2}kR\right)\left(1+\frac{1}{kR}+\frac{1}{k^2R^2}\right)$$ for some constant $C_0$ that is independent of $R$ and $k$. To show the stability result~(\ref{stab2}) recall from the second stability estimate  (3.5) of \cite{HMP11}
\begin{eqnarray*}
|u|_{2,\widetilde{\Omega}}&\leq& C_0\left(1+kR\right)\|\widetilde{f}\|_{\!_{0,\widetilde{\Omega}}},
\end{eqnarray*}
where $C$ is independent of $k$ and $u$.
Combining with the results above gives
\begin{eqnarray*}
|u|_{{2,\Omega}} &\leq& C_2\|f\|_{\!_{0,\Omega}}
\end{eqnarray*}
where $$C_2:=C_0(1+kR)C_1.\;\;\;\;\ep$$

\subsection{Some properties of Hankel functions}
We point out some properties of the DtN map that will be needed to derive error estimates. The result in Lemma \ref{hankel} can be found, for example, in Lemma 3.3 \cite{melenksauter}. 
\begin{lemma}\label{hankel}
For any $w\in H^\hf(\Gamma_{\!_R})$, it holds that
\begin{equation}
\begin{aligned}
&-\Imag\int_{\Gamma_{\!_R}}\DtN w\overline{w}\;ds \geq 0.
\end{aligned}
\end{equation}
\end{lemma}

The following lemma is identical to Lemma 3.2 of~\cite{Hsiao}, but in contrast to the result in~\cite{Hsiao}, we give the dependence of the constant on the wavenumber.

\begin{lemma}
For any $m\in\mathbb{Z}\backslash\{0\}$ and $kR>0$, it holds that 
\begin{equation}\label{hankelpositive}
\begin{aligned}
&\frac{1}{(1+m^2)^\hf}\left|\frac{H^{(2)'}_m(kR)}{H^{(2)}_m(kR)}\right|\leq \left|\frac{1}{m}\frac{H^{(2)'}_m(kR)}{H^{(2)}_m(kR)}\right|\leq \frac{1}{|m|}+\frac{1}{kR}.
\end{aligned}
\end{equation}
For $m=0$, we have that 
\begin{equation}
\begin{aligned}
&\left|\frac{H^{(2)'}_0(kR)}{H^{(2)}_0(kR)}\right|\leq C\left(1+\frac{1}{kR}\right),
\end{aligned}
\end{equation}
for some $C$ independent of $kR$.
\end{lemma}

{\bf Proof:}\\
We make use of the inequality (2.7) in \cite{wildemonk}
\begin{equation*}
\begin{aligned}
\left|H^{(2)}_\nu(\rho)\right|^2\left(\nu^2-\rho^2\right)+\rho^2\left|H^{(2)'}_\nu(\rho)\right|^2-\frac{4\rho}{\pi}\leq 0
\end{aligned}
\end{equation*}
and inequality (2.4) in \cite{wildemonk}
\begin{equation}
\begin{aligned}
\rho\left|H^{(2)}_\nu(\rho)\right|^2\geq\frac{2}{\pi}
\end{aligned}
\end{equation}
for any $\nu$ and $\rho\in\mathbb{R}$ such that $\nu\geq\hf$ and $\rho>0$. To obtain inequality~(\ref{hankelpositive}) for $m\geq 1$, choose $\nu = m$, $\rho=kR$. When $m\leq -1$, the result follows from the identity $$H^{(2)}_m(\rho)=(-1)^{|m|} H^{(2)}_{-m}(\rho).$$

To prove the inequality for $m=0$, recall the asymptotic relations of Hankel functions with small argument as $\rho\rightarrow 0$ (see e.g. 10.2.6, 10.7.2, 10.7.7 of \cite{NISTHandbook})
\begin{equation}\label{asymptotics0}
\begin{aligned}
H_\nu^{(2)}(\rho)\sim\begin{cases}
                   -(2i/\pi)\ln\rho\;\;\text{if $\nu=0$}\\
                   (2i/\pi)\Gamma(1)z^{-1}\;\;\text{if $\nu=1$},                    
                    \end{cases}
\end{aligned}
\end{equation}
and for large argument as $\rho\rightarrow\infty$
\begin{equation}\label{asymptotics1}
\begin{aligned}
H_\nu^{(2)}(\rho)\sim\begin{cases}
                   \sqrt{2/(\pi\rho)}e^{-i(\rho-\pi/4)}\;\;\text{if $\nu=0$}\\
                   \sqrt{2/(\pi\rho)}e^{-i(\rho-3\pi/4)}\;\;\text{if $\nu=1$}.                    
                    \end{cases}
\end{aligned}
\end{equation}
It follows from (\ref{asymptotics0}) and (\ref{asymptotics1}) that 
\begin{equation}
\begin{aligned}
\left|\frac{H^{(2)'}_0(\rho)}{H^{(2)}_0(\rho)}\right|=\left|\frac{H^{(2)}_1(\rho)}{H^{(2)}_0(\rho)}\right|\sim\begin{cases}
              C\frac{1}{\rho|\ln\rho|}\;\;\text{as $\rho\rightarrow 0$}\\
              1\hspace{1cm}\text{as $\rho\rightarrow\infty$}.
              \end{cases}
\end{aligned}
\end{equation}
Since $1/(\rho|\ln\rho|)<1/\rho$ for small enough $\rho$, we have by continuity of the Hankel functions in any bounded interval
\begin{equation}\label{hankelzero}
\begin{aligned}
\left|\frac{H^{(2)'}_0(\rho)}{H^{(2)}_0(\rho)}\right|\leq C\left(1+\frac{1}{\rho}\right).
\end{aligned}
\end{equation}
when $\rho>0$.
\scalebox{1.0}{$\square$}

\subsection{Truncated Boundary Value Problem}
In practical computations, one needs to truncate the infinite series of the DtN operator to obtain an approximate mapping written as a finite sum
\begin{equation}\label{tdtnmap}
\begin{aligned}
\DtN_{\!_N}w&=\sum_{|m|\leq N} k\frac{H^{(2)'}_m(kR)}{H^{(2)}_m(kR)}w_m e^{im\theta}
\end{aligned}
\end{equation}
for all $w\in H^\hf(\Gamma_{\!_R})$. The boundary value problem (\ref{dtnbvp}) is replaced by the following modified boundary value problem with a truncated DtN map: Find $u^{\!_N}\in H^1(\Omega)$ such that 
\begin{equation}\label{tdtnbvp}
\begin{aligned}
&\Delta\un + k^2\un = 0\;\;\mbox{in}\;\;\Omega\\
&\un = g\;\;\mbox{on}\;\;\Gamma_{\!_D}\\
&\diff{\un}{\bn}-\ThN\un=0\;\;\mbox{on}\;\;\Gamma_{\!_R}.
\end{aligned}
\end{equation}
In Theorem 4.5 of~\cite{Hsiao} it is shown that the truncated exterior \emph{Neumann} problem is well-posed for all $N$ sufficiently large. Following the same arguments, we can prove the following theorem for the truncated exterior \emph{Dirichlet} problem~(\ref{tdtnbvp}). 
\begin{theorem}\label{hpnx}
There exists an integer $N_0\geq 0$ depending on $k$ such that for any $g\in H^\hf(\Gamma_{\!_D})$ the truncated Dirichlet boundary value problem \emph{(\ref{tdtnbvp})} has a unique solution, $u^{\!_N}\in H^1(\Omega)$ for $N\geq N_0$.
\end{theorem}

\subsection{Stability of the truncated adjoint problem}
For later use in the derivation of $L^2$ norm estimates for the discretization error, we make use of duality arguments, which in turn depend on the following adjoint problem to~(\ref{tdtnbvp}):
\begin{equation}\label{truncatedadjoint}
\begin{aligned}
&\Delta z^{\!_N}+k^2 z^{\!_N}= f\;\;\mbox{in}\;\;\Omega,\\
&z^{\!_N} = 0\;\;\mbox{on}\;\;\Gamma_{\!_D},\\
&\diff{z^{\!_N}}{\bn}-\DtN_{\!_N}^{\star} z^{\!_N} = 0,
\end{aligned}
\end{equation}
where $f\in L^2(\Omega)$. The following lemma gives a stability constant $C_3$ of the solution $z^{\!_N}$ in the $H^2$ norm with respect to the $L^2$ norm of the data $f$. We note in particular that $C_3$ is independent of $N$. This first inequality in Lemma~\ref{reg3} below is proved in Lemma 4.2 of \cite{kmdtn}. The second inequality is proved in \cite{HMP15} before equation (29).
\begin{lemma}\label{reg3}
Assume $N\geq N_0$ is large enough so that the truncated boundary value problem~(\ref{truncatedadjoint}) has a unique solution $z^{\!_N}\in H^2(\Omega)$. Then there exist positive numbers $C_3$ and $C_4$ independent of $N$, $z^{\!_N}$ and $f$, but depending on $k$ and $R$ such that 
\begin{eqnarray}\label{stab4}
\|z^{\!_N}\|_{2,\Omega}&\leq& C_3\|f\|_{{0,\Omega}}\\
\|z^{\!_N}\|_{L^\infty(\Omega)}&\leq& C_4\|f\|_{0,\Omega}.
\end{eqnarray}
\end{lemma}

\section{The Plane Wave Discontinuous Galerkin Method}
\label{pwdg.sec}
Let $\mathscr{T}_h$ denote a finite element partition of $\Omega$ into elements $\{K\}$. We shall assume that all the elements $K\in\mathscr{T}_h$ are generalized triangles. A generalized triangle will be a true triangle in the interior of $\Omega$ but may have one curvilinear edge if the triangle has a face on $\Gamma_{\!_D}$ or $\Gamma_{\!_R}$. More general elements (e.g. quadrilaterals, pentagons, etc) are possible. The parameter $h$ represents the diameter of the largest element in $\mathscr{T}_h$, so that $h = \displaystyle\max_{K\in\mathscr{T}_h} h_{\!_K}$ where $h_{\!_K}$ is the diameter of the smallest circumscribed circle containing $K$. Denote by $\mathcal{E}$ the \emph{mesh skeleton}, i.e. the set of all edges of the mesh, $\cEI$ the set of interior edges, $\cED$ the set of edges on the boundary of the scatterer $\Gamma_{\!_D}$ and $\cEB$ the set of edges on the artificial boundary $\Gamma_{\!_R}$.

As is standard in DG methods, we introduce the \emph{jumps} and \emph{averages } as follows. Let $K^+,K^-\in\mathscr{T}_h$ be two elements sharing a common edge $e$. Suppose $\bn^+, \bn^-$ are the outward pointing unit normal vectors on the boundaries $\partial K^+$ and $\partial K^-$ respectively. Let $v:\Omega\rightarrow \mathbb{C}$ be a smooth scalar valued piecewise defined function,  and $\bsig:\Omega\rightarrow \mathbb{C}^2$ a smooth vector valued piecewise defined function.
Let  $\bx$ be a point on $e$. Define 
\begin{eqnarray*}
v^+(\bx) &:=& \lim_{{\bf y}\rightarrow \bx \atop{\bf y}\in K^+} v({\bf y}).
\end{eqnarray*}
The definitions of $v^-$, $\bsig^+$ and $\bsig^-$ are similar. The jumps are defined as 
\begin{eqnarray}\label{jumps}
&&\jmp{v}:= v^+\bn^+ + v^-\bn^-,\;\;\;\jmp{\bsig}:= \bsig^+\cdot\bn^+ + \bsig^-\cdot\bn^-,
\end{eqnarray}
and the averages are defined as 
\begin{eqnarray}\label{averages}
&&\avg{v}:=\hf\left(v^+ + v^-\right),\;\;\;\avg{\bsig}:=\hf\left(\bsig^+ +\bsig^-\right).
\end{eqnarray}
Suppose $u\in H^{2}(\Omega)$ is the exact solution of the homogeneous Helmholtz equation. In each element $K\in\mathscr{T}_h$, the weak formulation is:
\begin{eqnarray}\label{weakderv1}
\int_K\left(\nabla u\cdot\overline{\nabla v}-k^2u\overline{v}\right)\;d\bx -\int_{\partial K}\nabla u\cdot\bn_{\!_K}\overline{v}\;ds&=&0,
\end{eqnarray}
where $v$ is assumed piecewise smooth, and $\bn_{\!_K}$ is the outward pointing unit normal vector on $\partial K$. This smoothness assumption allows us to take traces of $v$ and $\nabla v$ on $\partial K$. Integrating equation~(\ref{weakderv1}) by parts once more leads to
\begin{eqnarray}\label{weakderv2}
\int_K u\overline{\left(-\Delta v-k^2v\right)}\;d\bx+\int_{\partial K}u\overline{\nabla v\cdot\bn_{\!_K}}\;ds-\int_{\partial K}\nabla u\cdot\bn_{\!_K}\overline{v}\;ds&=&0.
\end{eqnarray}
To proceed, we suppose the test function $v$ belongs in the Trefftz space $T(\mathscr{T}_h)$ defined as follows: Let $H^s(\mathscr{T}_h)$ be the broken Sobolev space on the mesh 
$$ H^s(\mathscr{T}_h):=\left\lbrace v\in L^2(\Omega): v|_{K}\in H^s(K)\;\;\forall\;K\in \mathscr{T}_h\right\rbrace.$$
Then the {Trefftz} space $T(\mathscr{T}_h)$ is
$$T(\mathscr{T}_h):=\left\lbrace v\in L^2(\Omega): v\in H^{2}(\mathscr{T}_h)\;\mbox{and}\;\Delta v+ k^2v=0\;\mbox{in\;each}\;K\in \mathscr{T}_h \right\rbrace.$$
Because $\Delta v+k^2v=0$ in $K$, equation~(\ref{weakderv2}) reduces to
\begin{eqnarray}\label{weakderv3}
\int_{\partial K}u\overline{\nabla v\cdot\bn_{\!_K}}\;ds-\int_{\partial K}\nabla u\cdot\bn_{\!_K}\overline{v}\;ds&=&0.
\end{eqnarray}
The problem now is to find an approximation of $u$ in a \emph{finite dimensional} Trefftz subspace of $T(\mathscr{T}_h)$. Define the finite dimensional \emph{local} solution space $V_{p_{\!_K}}(K)$ of dimension $p_{\!_K}\geq 1$ on each element $K\in\mathscr{T}_h$:
$$V_{{p_{\!_K}}}(K):=\left\lbrace w_h\in H^2(K):\Delta w_h +k^2w_h=0\;\;\mbox{in}\;K \right\rbrace$$
and the global solution space
$$V_h(\mathscr{T}_h):=\left\lbrace v\in L^2(\Omega): v|_{\!_K}\in V_{p_{\!_K}}(K)\;\;\mbox{in each $K\in\mathscr{T}_h$} \right\rbrace  $$
where the local dimension $p_{\!_K}$ can change from element to element.

Suppose that in each element $K\in\mathscr{T}_h$, $u_h$ is the unknown approximation of $u$ in the local solution space $V_{p_{\!_K}}(K)$ and $ik\bsig_h:=\nabla u_h$ is the flux. Then, on $\partial K$, we write 
\begin{eqnarray}\label{weakderv4}
\int_{\partial K}u_h\overline{\nabla v\cdot\bn_{\!_K}}\;ds-\int_{\partial K}ik\bsig_h\cdot\bn_{\!_K}\overline{v}\;ds&=&0,
\end{eqnarray}
for all $v\in V_{p_{\!_K}}(K)$. At this stage, $u_h$ and $ik\bsig_h$ are multi-valued on an edge $e\subset\partial K_1\cap\partial K_2$, since the trace from $K_1$ could differ from that of $K_2$. To find a \emph{global} numerical solution in $V_h(\mathscr{T}_h)$, we need $u_h$ and $ik\bsig_h$ to be \emph{single valued} on each edge of the mesh. Thus, we introduce \emph{numerical fluxes} $
\widehat{u}_h$ and $\widehat{\bsig}_h$ that are single valued approximations of $u_h$ and $ik\bsig_h$ respectively on each edge.

In each element of the mesh $K\in\mathscr{T}_h$, it holds that
\begin{eqnarray}\label{deriv5}
\int_{\partial K} \widehat{u}_h\overline{\nabla v\cdot\bn_{\!_K}}\;ds - \int_{\partial K} ik\widehat{\bsig}_h\cdot\bn_{\!_K}\overline{v}\;ds &=& 0.
\end{eqnarray}

Integration by parts allows us to write a ``domain based" equation that is equivalent to (\ref{deriv5})
\begin{eqnarray}\label{volmderv}
&&\int_K \left(\nabla u_h\cdot\overline{\nabla v}-k^2u_h\overline{v}\right)\;d\bx +\int_{\partial K}(\widehat{u}_h - u_h)\overline{\nabla v\cdot\bn_{\!_K}} - \int_{\partial K} ik\widehat{\bsig}_h\cdot\bn_{\!_K}\overline{v}\;ds=0.\nonumber\\
\end{eqnarray}
The form (\ref{volmderv}) is used to prove coercivity properties of the PWDG method, while the skeleton-based form (\ref{deriv5}) is used to program the method.

We now define the PWDG fluxes. The definition of the fluxes on interior edges and edges on the scatterer are taken to be those of standard PWDG methods in \cite{HMP13} and \cite{HMP11} as given in~(\ref{intfluxes}) and (\ref{dirfluxes}). But the fluxes on the artificial boundary $\Gamma_R$ are new. 
Denote by $\nabla_h$ the elementwise application of the gradient operator. Following \cite{HMP13},
\begin{equation}\label{intfluxes}
\left. \begin{aligned}
\widehat{u}_h &= \avg{u_h} - \frac{\beta}{ik}\jmp{\nabla_h u_h},\\
ik\widehat{\bsig_h} &= \avg{\nabla_h u_h} -\alpha ik \jmp{u_h}
\end{aligned}
\right\}
\qquad \text{on interior edges $\cEI$.}
\end{equation}
On the boundary of the scatterer $\Gamma_{\!_D}$
\begin{equation}\label{dirfluxes}
\left. \begin{aligned}
\widehat{u}_h &= 0,\\
ik\widehat{\bsig_h} &= {\nabla_h u_h} -\alpha ik {u_h}
\end{aligned}
\right\}
\qquad \text{on Dirichlet edges $\cED$.}
\end{equation}
For edges on the artificial boundary $\cEB$, we propose
\begin{eqnarray}\label{dtntr}
\widehat{u}^{\!_N}_h &=& u^{\!_N}_h-\frac{\delta}{ik}\left( \nabla_h u^{\!_N}_h\cdot \bn - \mathcal{S}_{\!_N} u^{\!_N}_h\right), \\
ik\widehat{\bsig}^{\!_N}_h &=& \mathcal{S}_{\!_{\!_N}} u^{\!_N}_h\bn-\frac{\delta}{ik}\mathcal{S}^{\star}_{\!_{\!_N}}\left(\nabla_h u^{\!_N}_h - \mathcal{S}_{\!_{\!_N}} u^{\!_N}_h\bn\right),
\end{eqnarray}
where $\alpha, \beta, \delta>0$ are positive flux coefficients defined on the edges of the mesh, and $\mathcal{S}^{\star}_{\!_{\!_N}}$ is the $L^2(\Gamma_{\!_R})$-adjoint of $\mathcal{S}_{\!_{\!_N}}$, defined as $$\int_{\Gamma_{\!_R}} \mathcal{S}_{\!_{\!_N}}^{\star}v\overline{w}\;ds = \int_{\Gamma_{\!_R}}v\overline{\mathcal{S}_{\!_{\!_N}} w}\;ds.$$
 Substituting these fluxes into equation ($\ref{volmderv}$), and summing over all elements $K\in\mathscr{T}_h$, we obtain the following PWDG scheme: Find $u^{\!_N}_h\in PW(\mathscr{T}_h)$ such that for all $v_h\in PW(\mathscr{T}_h)$
\begin{eqnarray}
\mathscr{A}_{{\!_N}}(u^{\!_N}_h, v_h) &=& \mathscr{L}_h(v_h)\label{pwdg2}
\end{eqnarray}
where
\begin{eqnarray}
\mathscr{A}_{\!_{\!_N}}(u^{\!_N}_h,v_h)&:=& \int_{\Omega}\left( \nabla_h u^{\!_N}_h\cdot\nabla_h\overline{v_h}-k^2 u^{\!_N}_h \overline{v_h}\right)\;d{\bx}-\int_{\cEI}\jmp{u^{\!_N}_h}\cdot\overline{\avg{\nabla_h v_h}}\;ds\nonumber\\&&\textcolor{black}{-\int_{\cEB} \mathcal{S}_{\!_{\!_N}} u^{\!_N}_h\overline{v_h}}\;ds-\frac{1}{ik}\int_{\cEI}\beta \jmp{\nabla_h u^{\!_N}_h}\overline{\jmp{\nabla_h v_h}}\;ds\nonumber\\&&-\int_{\cEI}\avg{\nabla_h u^{\!_N}_h}\cdot\overline{\jmp{v_h}}\;ds +ik\int_{\cEI}\alpha\jmp{u^{\!_N}_h}\cdot\overline{\jmp{v_h}}\;ds\nonumber\\&& -\frac{1}{ik}\int_{\cEB}\delta\left(\nabla_h u^{\!_N}_h\cdot \bn - \mathcal{S}_{\!_{\!_N}}u^{\!_N}_h\right)\overline{\left(\nabla_h v_h\cdot \bn - \mathcal{S}_{\!_{\!_N}} v_h\right)}\;ds\nonumber\\
&&+{ik}\int_{\cED}\alpha  u^{\!_N}_h\overline{v_h}\;ds-\int_{\cED}\left(u^{\!_N}_h\overline{\nabla_h v_h\cdot\bn}+\nabla_h u^{\!_N}_h\cdot\bn \overline{v_h}\right)\;ds,\label{pwdg3}
\end{eqnarray}
and the right hand side is
\begin{eqnarray}
\mathscr{L}_h(v_h) &:=& -\int_{\cED} g\overline{\nabla_h v_h\cdot\bn}\;ds+{ik}\int_{\cED}\alpha g\overline{v_h}\;ds.\label{rhs1}
\end{eqnarray}
The hermitian sesquilinear form ($\ref{pwdg3}$) allows us to prove coercivity of the DtN-PWDG scheme. However to program the method, we can make the algorithm more efficient by exploiting the Trefftz property of $PW(\mathscr{T}_h)$ to write the sesquilinear form on the skeleton of the mesh. This avoids the need to integrate over elements in the mesh. Integrating by parts ($\ref{volmderv}$) and using the Trefftz property $\Delta v_h + k^2 v_h = 0$, the elemental equation (\ref{volmderv}) reduces to 
\begin{eqnarray}
&&\int_{\partial K} \widehat{u}^{\!_N}_h\overline{\nabla_h v_h\cdot\bn}\;ds - \int_{\partial K} ik\widehat{\boldsymbol \sigma}^{\!_N}_h\cdot\bn\overline{v}^{}_h\;ds = 0.
\end{eqnarray} Then substituting the numerical fluxes and summing over all elements of the mesh, we get
\begin{eqnarray}
\mathscr{A}_{\!_{\!_N}}(u^{\!_N}_h,v_h)&=& \int_{\cEI}\avg{u^{\!_N}_h}\overline{\jmp{\nabla_h v_h}}\;ds-\int_{\cEI}\avg{\nabla_h u^{\!_N}_h}\cdot\overline{\jmp{v_h}}\;ds-\int_{\cEB}\mathcal{S}_{\!_{\!_N}} u^{\!_N}_h\overline{v_h}\;ds\nonumber \\
&& +\int_{\cEB}u^{\!_N}_h\overline{\nabla_h v_h\cdot\bn}\;ds-\frac{1}{ik}\int_{\cEI}\beta \jmp{\nabla_h u^{\!_N}_h}\overline{\jmp{\nabla_h v_h}}\;ds\nonumber\\
&& +ik\int_{\cEI}\alpha\jmp{u^{\!_N}_h}\cdot\overline{\jmp{v_h}}\;ds
-\int_{\cED} \nabla_h u^{\!_N}_h \cdot\bn\overline{v_h}\;ds+{ik}\int_{\cED}\alpha  u^{\!_N}_h\overline{ v_h}\;ds \nonumber \\ && -\frac{1}{ik}\int_{\cEB}\delta\left(\nabla_h u^{\!_N}_h\cdot\bn - \mathcal{S}_{\!_{\!_N}}u^{\!_N}_h\right)\overline{\left(\nabla_h v_h\cdot \bn - \mathcal{S}_{\!_{\!_N}}v_h\right)}\;ds
 \label{pwdg4}
\end{eqnarray}
Our DtN-PWDG MATLAB code is based on the sesquilinear form (\ref{pwdg4}).

For error estimates, it is useful to derive an equivalent form of ~(\ref{pwdg3}). Applying a DG magic formula (Lemma 6.1, \cite{MelenkEsterhazy12}) to ~(\ref{pwdg3})
we get
\begin{eqnarray}
\mathscr{A}_{\!_{\!_N}}(u^{\!_N}_h,v_h)&=& \int_{\cEI}\jmp{\nabla_h u^{\!_N}_h}\overline{\avg{ v_h}}\;ds-\int_{\cEI}\jmp{u^{\!_N}_h}\cdot\overline{\avg{\nabla_h v_h}}\;ds-\int_{\cED} u^{\!_N}_h\overline{\nabla_h v_h\cdot\bn}\;ds \nonumber\\&&-\frac{1}{ik}\int_{\cEI}\beta \jmp{\nabla_h u^{\!_N}_h}\overline{\jmp{\nabla_h v_h}}\;ds +ik\int_{\cEI}\alpha\jmp{u^{\!_N}_h}\cdot\overline{\jmp{v_h}}\;ds\nonumber\\&&
+ {ik}\int_{\cED}\alpha u^{\!_N}_h \overline{v_h}\;ds\textcolor{black}{+\int_{\cEB} \left(\nabla_h u^{\!_N}_h\cdot\bn-\mathcal{S}_{\!_{\!_N}} u^{\!_N}_h\right)\overline{v_h}}\;ds\nonumber \\ && -\frac{1}{ik}\int_{\cEB}\delta\left(\nabla_h u^{\!_N}_h\cdot\bn - \mathcal{S}_{\!_{\!_N}}u^{\!_N}_h\right)\overline{\left(\nabla_h v_h\cdot\bn - \mathcal{S}_{\!_{\!_N}} v_h\right)}\;ds.\label{pwdg5}
\end{eqnarray}

 \begin{proposition}\label{pwdgcons}
The \emph{DtN-PWDG} method is consistent.
\end{proposition}
{\bf Proof.}
If $u^{\!_N}$ is the exact solution of the truncated boundary value problem $(\ref{tdtnbvp})$, then under the assumptions on the geometry of the scatterer, $u^{\!_N}\in H^{2}(\Omega)$, thus on any interior edge $e$, $\jmp{u^{\!_N}}={\bf 0}$ and $\jmp{\nabla_h u^{\!_N}}=0$ on $\cEI$, $\nabla_h u^{\!_N}\cdot \bn= {\mathcal{S}}_{\!_N}u^{\!_N}$ on $\cEB$, and $u^{\!_N}=g$ on $\cED$. 
Therefore from {(\ref{pwdg5})}, for any $v\in PW(\mathscr{T}_h)$
\begin{eqnarray}
\mathscr{A}_{\!_N}(u^{\!_N},v) &=& -\int_{\cED} g\overline{\nabla_h v\cdot\bn}\;ds +ik\int_{\cED}\alpha g\overline{v}\;ds\nonumber\\
&=& \mathscr{L}_h(v).\;\;\;\scalebox{1.0}{$\square$}
\end{eqnarray} 
\begin{proposition}
Provided $N\geq N_0$, the mesh-dependent functional
\begin{eqnarray}
\|v\|_{\!_{DG}}&:=& \sqrt{\Imag{\mathscr{A}_{\!_{\!_N}}(v,v)}} \label{dgnorm}
\end{eqnarray}
defines a norm on $T(\mathscr{T}_h)$. Moreover, setting 
\begin{eqnarray}
\|v\|^2_{\!_{\!_{DG^+}}} &:=& \|v\|^2_{\!_{\!_{DG}}}+k\|\beta^{-\frac{1}{2}}\avg{v}\|^2_{0,\cEI}+k^{-1}\|\alpha^{-\frac{1}{2}}\avg{\nabla_h v}\|^2_{0,\cEI} \nonumber\\
&& +k^{-1}\|\alpha^{-\frac{1}{2}}{\nabla_h v\cdot\bn}\|^2_{0,\cED}+k\|\delta^{-\frac{1}{2}}v\|^2_{0,\cEB} 
\end{eqnarray}
we have 
\begin{eqnarray}\label{cont}
\mathscr{A}_{\!_{\!_N}}(v,w)\leq 2\|v\|_{\!_{\!_{DG}}}\|w\|_{\!_{\!_{DG^+}}}
\end{eqnarray}

\end{proposition}

{\bf Proof.}
Taking the imaginary part of (\ref{pwdg3}), we have
\begin{eqnarray}
\Imag \mathscr{A}_{\!_{\!_N}}(v,v)&=& k^{-1}\|\beta^{\frac{1}{2}}\jmp{\nabla_h v}\|^2_{0,\cEI}+k\|\alpha^{\frac{1}{2}}\jmp{v}\|^2_{0,\cEI}+k\|\alpha^{\frac{1}{2}}{ v}\|^2_{0,\cED}-\Imag\int_{\cEB}\mathcal{S}_{\!_{\!_N}} v\overline{v}\;ds\nonumber\\
&& + k^{-1}\| \delta^{1/2} (\nabla_h v\cdot\bn-\mathcal{S}_{\!_{\!_N}}v)\|^2_{0,\cEB}\nonumber\\
&=& \|v\|^2_{\!_{DG}}.
\end{eqnarray}
From Lemma~\ref{hankel}, we recall that taking only partial sums, 
\begin{eqnarray*}
&&-\Imag\int_{\cEB}\mathcal{S}_{\!_{\!_N}} v\overline{v}\;ds=\sum_{|m|\leq N} \frac{4|v_m|^2}{|H^{(2)}_m(kR)|^2}>0
\end{eqnarray*}
If $\Imag \mathscr{A}_{\!_{\!_N}}(v,v) = 0$, then $v\in H^{2}(\Omega)$ satisfies the Helmholtz equation $\Delta v+k^2v=0$ in $\Omega$, with $v=0$ on $\Gamma_{\!_D}$, and $\nabla v\cdot\bn -\mathcal{S}_{\!_{\!_N}}v=0$ on $\Gamma_{\!_R}$ (our new choice of flux is key to asserting this last equality). By Theorem~\ref{hpnx}, this problem has only the trivial solution $v=0$ provided $N\geq N_0$ is large enough.

To prove (\ref{cont}), we apply the Cauchy Schwarz inequality repeatedly to {(\ref{pwdg5})}.
\scalebox{1.0}{$\square$}

{\bf Remark:} The assumption that $0<\delta\leq\hf$ required to prove continuity of the sesquilinear form for the PWDG with impedance boundary conditions as in \cite{HMP13, HMP11} is no longer necessary. It is sufficient that $\delta>0$ in the DtN-PWDG scheme. This is because of our new choice of boundary fluxes.

\begin{proposition}
Provided $N\geq N_0$, the discrete problem \emph{(\ref{pwdg2})} has a unique solution $u^{\!_N}_h\in PW(\mathscr{T}_h)$.
\end{proposition}
{\bf Proof.}
Assume $\mathscr{A}_{\!_N}(u^{\!_N}_{h},v) = 0$ for all $v\in PW(\mathscr{T}_h)$. Then in particular $\mathscr{A}_{\!_N}(u^{\!_N}_{h},u^{\!_N}_{h}) =0$ and so $\Imag \mathscr{A}_{\!_N}(u^{\!_N}_h,u^{\!_N}_h)=0$. Then $\|u^{\!_N}_h\|_{\!_{DG}}=0$ which implies $u^{\!_N}_h=0$ since $\|\cdot\|_{\!_{DG}}$ is a norm on $PW(\mathscr{T}_h)$.
\ep

\section{Error Estimates}
\label{error.sec}
\subsection{A mesh-dependent error estimate}

 We state an error estimate in the mesh dependent $DG$ norm.

\begin{proposition}\label{eest1}
Assume $N\geq N_0$. Let $u^{\!_N}$ be the unique solution of the truncated boundary value problem~\emph{(\ref{tdtnbvp})}, and $u^{\!_N}_h\in PW(\mathscr{T}_h)$ the unique solution of the discrete problem~\emph{(\ref{pwdg2})}. Then
\begin{eqnarray}
\|u^{\!_N} - u^{\!_N}_h\|_{\!_{DG}}&\leq& 2\inf_{w^{}_h\in PW(\mathscr{T}_h)}\|u^{\!_N} - w^{}_h\|_{\!_{DG^+}}.
\end{eqnarray}
\end{proposition}

{\bf Remark:} Using Moiola's estimates for approximation by plane waves, this result can be used to provide order estimates (see Corollary \ref{coro: kapita2} later in this paper).

{\bf Proof}: Let $w^{\!_N}\in\pwdg{\mathscr{T}_h}$ be arbitrary. By Proposition~\ref{pwdgcons}, definition (\ref{dgnorm}) and the inequality~(\ref{cont}) we have
\begin{eqnarray*}
\|u^{\!_N}-u_h^{\!_N}\|^2_{\!_{DG}} &=& \Imag \mathscr{A}_{\!_N}(u^{\!_N}-u^{\!_N}_h,u^{\!_N}-u^{\!_N}_h)\\
&\leq& |\mathscr{A}_{\!_N}(u^{\!_N}-u^{\!_N}_h,u^{\!_N}-u^{\!_N}_h)|\\
&=& |\mathscr{A}_{\!_N}(u^{\!_N}-u^{\!_N}_h,u^{\!_N}-w^{}_h)|\\
&\leq& 2\|u^{\!_N}-u^{\!_N}_h\|_{\!_{DG}}\|u^{\!_N}-w^{}_h\|_{\!_{DG^+}}.\;\;\;\;\ep
\end{eqnarray*}

\subsection{Error estimates in the $L^2$ norm}

Let $e^{\!_N}_h=u-u^{\!_N}_h$ be the error of the DtN-PWDG method, where $u$ is the exact solution, and $u^{\!_N}_h$ the computed solution. Let $u^{\!_N}$ be the solution of the truncated boundary value problem~(\ref{tdtnbvp}). Then, by the triangle inequality
\begin{eqnarray}
\|u-u^{\!_N}_h\|_{L^2(\Omega)}&\leq& \|u-u^{\!_N}\|_{L^2(\Omega)} + \|u^{\!_N}-u^{\!_N}_h\|_{L^2(\Omega)}.
\end{eqnarray}
The term $\|u-u^{\!_N}\|_{L^2(\Omega)}$ is the \emph{truncation error} introduced by truncating the DtN map, while the term $\|u^{\!_N}-u^{\!_N}_h\|_{L^2(\Omega)}$ is the \emph{discretization error} of the DtN-PWDG method.

\subsubsection{Estimation of $\|u-u^{\!_N}\|_{L^2(\Omega)}$}\label{sec:truncationerror}

 Error estimates for the truncation error of the Helmholtz problem (\ref{dtnbvp}) with DtN boundary conditions were proved by D. Koyama (see \cite{kdtn, kmdtn}), where it is assumed that $\mbox{supp}(f)\subset \Omega^{\prime}$, $\Omega^{\prime}\subset B_a$ where $B_{a}$ is a ball of radius $a\leq R$, such that $B_a\subset B_R$.
 
The main result of this section is the following (see Proposition 4.1 of \cite{kmdtn}).
\begin{theorem}\label{trthm}
Assume $N\geq N_0$ and that the exact solution $u$ of problem (\ref{helmholtzbc}) belongs to $H^{m}(\Omega)$. Assuming that $\mbox{supp}(f)\subset \Omega^{\prime}$, $\Omega^{\prime}\subset B_a$ where $B_{a}$ is a ball of radius $a\leq R$, such that $B_a\subset B_R$, there exists a positive number $C_4$ independent of $N$, $f$, $u$ and $u^{\!_N}$, but depending on $k,a,R$ and $\Omega$, such that for $m\in\mathbb{N}\cup\{0\}$ we have
$$\|u-u^{\!_N}\|_{m,\Omega^{\prime}} \leq  C_4 N^{-\hf-s+m}\left|\frac{H^{(2)}_{\!_N}(kR)}{H^{(2)}_{\!_N}(ka)}\right|\mathscr{R}_{\!_N}(u;s,a)$$
where $\Omega$ may be $\Omega^\prime$ and $$\mathscr{R}_{\!_N}(u;s,a):= \left(\sum_{|m|>N}|m|^{2s}|u_m(a)|^2\right)^\hf.$$
\end{theorem}

\subsubsection{Estimation of $\|u^{\!_N}-u^{\!_N}_h\|_{L^2(\Omega)}$} 
In this section, we study a priori error estimates for the discretization error
$\|u^{\!_N}-u^{\!_N}_h\|_{L^2(\Omega)}$ of the DtN-PWDG scheme.
For this section, we make the following assumptions on the mesh and flux parameters.

\emph{Assumption on the mesh}:
The mesh is \emph{shape regular} and \emph{quasi-uniform}. The method can be applied to more general meshes such as quadrilateral elements and locally refined meshes, but our focus in this chapter is the boundary condition on $\Gamma_{\!_R}$, so we choose simple quasi-uniform meshes. 

\emph{Assumption on the numerical fluxes}:
Since the mesh is quasi-uniform, we assume the numerical fluxes $\alpha,\beta,\delta$ are positive universal constants on the mesh.

Suppose $e\subset\partial K_1\cap\partial K_2$ is a common edge between triangles $K_1$ and $K_2$ (on boundary edges assume $K_2$ is empty). We recall the trace inequality (see \cite{brennerscott}, Theorem 1.6.6),
\begin{equation}\label{trace1}
\begin{aligned}
\|v\|^2_{0,e}&\leq C\sum_{j=1}^2\left(h^{-1}_{\!_{K_j}} \|v\|^2_{0,K_j}+h_{\!_{K_j}}|v|^2_{1,K_j} \right)
\end{aligned}
\end{equation}
with $C$ depending only on the shape regularity measure $\mu$. In addition, as $v\in H^{2}(\Omega)$, it holds that
\begin{equation}\label{trace2}
\begin{aligned}
\|\nabla v\|^2_{0,e} &\leq C\sum_{j=1}^2\left(h^{-1}_{\!_{K_j}}\|\nabla v\|^2_{0,K_j} +h_{\!_{K_j}}|v|^2_{2,K_j}\right)
\end{aligned}
\end{equation}
with $C$ depending only on the shape regularity parameter $\mu$.

Choosing $f =e^{\!_N}_h:= u^{\!_N}-u^{\!_N}_h$ in~(\ref{truncatedadjoint}), we have by the adjoint consistency of the DtN-PWDG scheme that
\begin{eqnarray}\label{adconsdtn}
\mathscr{A}_{\!_N}(w,z^{\!_N}) &=& \int_{\Omega} w\overline{e^{\!_N}_h}\;d\bx.
\end{eqnarray}
where $w\in T(\mathscr{T}_h)$ is any piecewise solution of the Helmholtz equation. Choosing $w=e^{\!_N}_h$ in~(\ref{adconsdtn}), the consistency of the DtN-PWDG method in Proposition~\ref{pwdgcons} implies that
\begin{eqnarray}
\|e^{\!_N}_h\|^2_{L^2(\Omega)}&=& \mathscr{A}_{\!_N}(e^{\!_N}_h,z^{\!_N})\nonumber\\
&=& \mathscr{A}_{\!_N}(e^{\!_N}_h,z^{\!_N}-z^{}_h)\label{l2error1}
\end{eqnarray}
for any arbitrary $z^{}_h\in PW(\mathscr{T}_h)$. To approximate the right hand side
of~(\ref{l2error1}) we follow the idea introduced in Lemma 5.3 of~\cite{kmw} and Lemma 3.10 of \cite{GHP09}: Let $z^c_h$ be the conforming piecewise linear finite element interpolant of $z\in H^{2}(\Omega)$. Then we can find a $z^{}_{h}\in PW(\mathscr{T}_h)$ that can approximate $z^c_h$. Adding and subtracting $z^c_h$, we have
\begin{eqnarray*}
\|e^{\!_N}_h\|^2_{L^2(\Omega)} &=&  \mathscr{A}_{\!_N}(e^{\!_N}_h,z^{\!_N}-z^c_h)+\mathscr{A}_{\!_N}(e^{\!_N}_h,z^c_h-z^{}_h)
\end{eqnarray*}
For convenience of notation, denote by $e_h^{\!_{N,c}}:=z^{\!_N}-z^c_h$ the conforming error.
It follows that 
\begin{eqnarray}\label{sesqerr1}
&&\mathscr{A}_{\!_N}(e^{\!_N}_h,e_h^{\!_{N,c}})= \int_{\cEI}\jmp{\nabla_h e^{\!_N}_h}\overline{\avg{e_h^{\!_{N,c}}}}\;ds-\int_{\cEI}\jmp{e^{\!_N}_h}\cdot\overline{\avg{\nabla_h (e_h^{\!_{N,c}})}}\;ds\\
&&-\frac{1}{ik}\int_{\cEI}\beta \jmp{\nabla_h e^{\!_N}_h}\overline{\jmp{\nabla_h (e_h^{\!_{N,c}})}}\;ds
+ik\int_{\cEI}\alpha\jmp{e^{\!_N}_h}\cdot\overline{\jmp{e_h^{\!_{N,c}}}}\;ds\nonumber\\&&
-\frac{1}{ik}\int_{\cEB}\delta\left(\nabla_h e^{\!_N}_h\cdot\bn - \DtN_{\!_N} e^{\!_N}_h\right)\overline{\left(\nabla_h e_h^{\!_{N,c}}\cdot\bn - \DtN_{\!_N} e_h^{\!_{N,c}}\right)}\;ds\nonumber \\
&&-\int_{\cED} e^{\!_N}_h\overline{\nabla_h e_h^{\!_{N,c}}\cdot\bn}\;ds \textcolor{black}{+\int_{\cEB} \left(\nabla_h e^{\!_N}_h\cdot\bn-\DtN_{\!_N} e^{\!_N}_h\right)\overline{e_h^{\!_{N,c}}}}\;ds+ {ik}\int_{\cED}\alpha e^{\!_N}_h \overline{e_h^{\!_{N,c}}}\;ds\nonumber. \label{bilinear}
\end{eqnarray}
As $e_h^{\!_{N,c}}$ is continuous, it follows that $\jmp{e_h^{\!_{N,c}}}=0$ on each interior edge, so that 
\begin{eqnarray}
ik\int_{\cEI}\alpha\jmp{e^{\!_N}_h}\cdot\overline{\jmp{e_h^{\!_{N,c}}}}\;ds &=& 0.
\end{eqnarray} 
Consider terms involving $e_h^{\!_{N,c}}$ in the sesquilinear form~(\ref{sesqerr1}):  
 \begin{eqnarray}
 I_1&=&\left|\int_{\cEI}\jmp{\nabla_h e^{\!_N}_h}\overline{\avg{e_h^{\!_{N,c}}}}\;ds\right|\nonumber\\
 &\leq& \sum_{e\in\cEI} k^{-\frac{1}{2}}\|\beta^{\frac{1}{2}}\jmp{\nabla_h e^{\!_N}_h}\|_{\!_{0,e}}\;k^{\frac{1}{2}}\|\beta^{-\frac{1}{2}}e_h^{\!_{N,c}}\|_{\!_{0,e}},\nonumber\\
 I_2 &=& \left|ik\int_{\cED}\alpha e^{\!_N}_h\cdot\bn\overline{e_h^{\!_{N,c}}}\;ds\right|\nonumber\\
 &\leq& \sum_{e\in\cED} k^{\frac{1}{2}}\|\alpha^{\frac{1}{2}}e^{\!_N}_h\|_{\!_{0,e}}\;k^{\frac{1}{2}}\|\alpha^{\frac{1}{2}}e_h^{\!_{N,c}}\|_{\!_{0,e}},\nonumber\\
 I_3 &=& \left|\int_{\cEB} \left(\nabla_h e^{\!_N}_h\cdot\bn-\DtN_{\!_N} e^{\!_N}_h\right)\overline{e_h^{\!_{N,c}}}\;ds\right|\nonumber\\
 &\leq& \sum_{e\in\cEB} k^{-\frac{1}{2}}\|\delta^{\frac{1}{2}}\left(\nabla_h e^{\!_N}_h\cdot\bn-{\DtN}_{\!_N}e^{\!_N}_h\right)\|_{\!_{0,e}}\;k^{\frac{1}{2}}\|\delta^{-\frac{1}{2}}e_h^{\!_{N,c}}\|_{\!_{0,e}}.\nonumber
 \end{eqnarray}
For any edge $e\in\cEI$, assume that $e\subset \partial K_1\cap \partial K_2$. Under the assumption that the numerical fluxes are constant, we have by the trace inequality (\ref{trace1}) that 
\begin{eqnarray}
\|\beta^{-\frac{1}{2}}e_h^{\!_{N,c}}\|_{\!_{0,e}}&\leq& C\sum_{\ell=1}^2 \left(\frac{1}{h_{K_\ell}^{\frac{1}{2}}}\|e_h^{\!_{N,c}}\|_{\!_{0,K_\ell}} + h_{K_\ell}^{\frac{1}{2}}\|\nabla_h e_h^{\!_{N,c}}\|_{\!_{0,K_\ell}}\right)\nonumber\\
&\leq& C\sum_{\ell=1}^2 h_{K_\ell}^{\htf}|z^{\!_N}|_{\!_{2,K_\ell}}.
\end{eqnarray}
 The other terms involving $\alpha,\delta$ are estimated in exactly the same way. Therefore,
 \begin{eqnarray}\label{iest}
 I_1+I_2+I_3 &\leq& C(\mu,\alpha,\beta,\delta,k,R)\|e^{\!_N}_h\|_{\!_{DG}}\sum_{K\in\mathscr{T}_h}h^{\htf}_{K}|z^{\!_N}|_{\!_{2,K}}
 \end{eqnarray}
 where $C$ depends on the wavenumber $k$, radius $R$, shape regularity parameter $\mu$ and the flux parameters  $\alpha,\beta$ and $\delta$, but is independent of $h$, $z^{\!_N}$ and $N$.
 
 Now consider terms involving $\nabla_h e_h^{\!_{N,c}}\cdot\bn$ in (\ref{bilinear}).
 \begin{eqnarray}
 J_1 &=& \left|\int_{\cEI}\jmp{e^{\!_N}_h}\cdot\overline{\avg{\nabla_h e_h^{\!_{N,c}}}}\;ds    \right|\nonumber\\
 &\leq& \sum_{e\in\cEI} k^{\frac{1}{2}}\|\alpha^{\frac{1}{2}}\jmp{e^{\!_N}_h}\|_{\!_{0,e}}\;k^{-\frac{1}{2}}\|\alpha^{-\frac{1}{2}}\avg{\nabla_h e_h^{\!_{N,c}}}\|_{\!_{0,e}},\nonumber\\
 J_2 &=& \left|-\frac{1}{ik}\int_{\cEI}\beta \jmp{\nabla_h e^{\!_N}_h}\overline{\jmp{\nabla_h e_h^{\!_{N,c}}}}\;ds \right|\nonumber\\
 &\leq& \sum_{e\in\cEI} k^{-\frac{1}{2}}\|\beta^{\frac{1}{2}}\nabla_h e^{\!_N}_h\|_{\!_{0,e}}\;k^{-\frac{1}{2}}\|\beta^{\frac{1}{2}}\nabla_h e_h^{\!_{N,c}}\|_{\!_{0,e}},\nonumber\\
 J_3 &=& \left|-\int_{\cED}e^{\!_N}_h\overline{{\nabla_h e_h^{\!_{N,c}}\cdot\bn}}\;ds \right|\nonumber\\
 &\leq& \sum_{e\in\cED} k^{\frac{1}{2}}\|\alpha^{\frac{1}{2}}e^{\!_N}_h\|_{\!_{0,e}}\;k^{-\frac{1}{2}}\|\alpha^{-\frac{1}{2}}\nabla_h e_h^{\!_{N,c}}\|_{\!_{0,e}}\nonumber,\\
 J_4 &=& \left| -\frac{1}{ik}\int_{\Gamma_{\!_{R}}}\delta\left(\nabla_h e^{\!_N}_h\cdot\bn - \DtN_{\!_N} e^{\!_N}_h\right)\overline{\left(\nabla_h e_h^{\!_{N,c}}\cdot\bn - \DtN_{\!_N} e_h^{\!_{N,c}}\right)}\;ds \right|\nonumber\\
 &\leq&  k^{-\frac{1}{2}}\|\delta^{\frac{1}{2}}(\nabla_h e^{\!_N}_h\cdot\bn-{\DtN}_{\!_N}e^{\!_N}_h)\|_{\!_{0,\Gamma_{\!_R}}}\;k^{-\frac{1}{2}}\|\delta^{\frac{1}{2}}(\nabla_h e_h^{\!_{N,c}}\cdot\bn-{\DtN}_{\!_N}e_h^{\!_{N,c}})\|_{\!_{0,\Gamma_{\!_R}}}.\nonumber 
 \end{eqnarray} 
 We start by approximating $J_1$. By the trace inequality (\ref{trace2}), we have on an edge $e\subset \partial K_1\cap \partial K_2$,
 \begin{eqnarray}
 \|\alpha^{-\frac{1}{2}}\avg{\nabla_h e_h^{\!_{N,c}}}\|_{\!_{0,e}}&\leq& C\sum_{\ell=1}^2\left(\frac{1}{h_{K_\ell}^{\frac{1}{2}}}\|\nabla_h e_h^{\!_{N,c}}\|_{\!_{0,K_\ell}}+h_{K_\ell}^{\hf}|e_h^{\!_{N,c}}|_{{2,K_\ell}}\right)\nonumber\\
 &\leq& C\sum_{\ell=1}^2 h_{K_\ell}^{\hf}|z^{\!_N}|_{{2,K_\ell}}.\nonumber
 \end{eqnarray}  
 The same argument holds for $J_2$ and $J_3$. The term $J_4$ that involves the DtN map can be estimated via the triangle inequality 
 \begin{eqnarray}
 &&k^{-\frac{1}{2}}\|\delta^{\frac{1}{2}}(\nabla_h e_h^{\!_{N,c}}\cdot\bn-{\DtN}_{\!_N} e_h^{\!_{N,c}})\|_{\!_{0,\Gamma_{\!_R}}}\nonumber\\&& \;\;\;\leq k^{-\frac{1}{2}}\|\delta^{\frac{1}{2}}\nabla_h e_h^{\!_{N,c}}\|_{\!_{0,\Gamma_{\!_R}}} + k^{-\frac{1}{2}}\|\delta^{\frac{1}{2}}{\DtN}_{\!_N}e_h^{\!_{N,c}}\|_{\!_{0,\Gamma_{\!_R}}}.
 \label{gammart}
 \end{eqnarray}
 The first term in (\ref{gammart}) is estimated as before. Now we estimate the term with the DtN map.  
 For $m\in\mathbb{Z}$, define $m_{\!_0}$ as follows 
 \begin{eqnarray}\label{mdef}
 m_{\!_0}&:=& \begin{cases} m &\mbox{if}\; |m|\neq 0 \\ 
1 & \mbox{if}\; m=0.\end{cases}
\end{eqnarray}
Let $$\mathscr{T}^{\!_R}_{h}:=\left\lbrace K\in\mathscr{T}_h: \mbox{length}\left(\partial K\cap \Gamma_R\right)>0 \right\rbrace $$ be the set of all elements with an edge on $\Gamma_{\!_R}$. Then the term in $J_4$ with the DtN map is estimated as
 \begin{eqnarray}
 \|\delta^{\frac{1}{2}}{\DtN}_{\!_N}(z^{\!_N}-z^c_h)\|_{\!_{0,\Gamma_{\!_R}}} &=& 2\pi \delta^\hf \left(\sum_{|m|\leq N} |m_{\!_0}|^2\frac{k^2R^2}{|m_{\!_0}|^2}\left|\frac{H^{(2)'}_{m}(kR)}{H^{(2)}_{m}(kR)}\right|^2|(e_h^{\!_{N,c}})_{m}|^2\right)^{\frac{1}{2}}\nonumber\\
 &\leq& C(\delta,k,R)\;N \|e_h^{\!_{N,c}}\|_{0,\Gamma_{\!_R}}\nonumber\\
 &\leq& C(\delta,k,R)\;N\sum_{e\in\cEB} \|e_h^{\!_{N,c}}\|_{0,e}\nonumber\\
 &\leq& C(\mu,\delta,k,R)\;N\sum_{K\in\mathscr{T}^{\!_R}_h}h_K^{\hf}|z^{\!_N}|_{{2,K}}.
 \end{eqnarray} 
Hence, summarizing we have the estimate  
 \begin{eqnarray}\label{jest}
 J_1+J_2+J_3+J_4 &\leq& Ch^\hf\left(N+1\right)\|e^{\!_N}_h\|_{\!_{DG}} \displaystyle\sum_{K\in\mathscr{T}_h}h_K^{\hf}|z^{\!_N}|_{{2,K}}\nonumber\\
 \end{eqnarray}
 
 Combining (\ref{iest}), (\ref{jest}) and (\ref{reg3}), we arrive at 
 {}{\begin{eqnarray}\label{bilinear1}
\mathscr{A}_{\!_N}(e^{\!_N}_h,z^{\!_N}-z^c_h)&\leq& Ch^\hf(N+1)|z^{\!_N}|_{{2,\Omega}} \|e^{\!_N}_h\|_{\!_{DG}}\nonumber\\
 &\leq&  Ch^\hf(N+1)\|e^{\!_N}_h\|_{\!_{DG}}\|u-u^{\!_N}_h\|_{0,\Omega}.
 \end{eqnarray} }
 where we have used the quasi-uniformity of the mesh and the stability estimate~(\ref{stab4}).
 
 To estimate the term $\mathscr{A}_{\!_N}(e^{\!_N}_h,z^c_h-z^{}_h)$, we use results from Lemma 5.4 of \cite{kmw}. By similar arguments used in the estimation of $\mathscr{A}_{\!_N}(e^{\!_N}_h,z^{\!_N} -z^{c}_h)$, and using the second inequality of Lemma \ref{reg3} to bound $\|z^{\!_N}\|_{L^\infty(\Omega)}$, we get
{}{ \begin{eqnarray}\label{bilinear2}
\mathscr{A}_{\!_N}(e^{\!_N}_h,z^c_h-z^{\!_N}) &\leq& Ch^{\htf}(N+1)\|e^{\!_N}_h\|_{\!_{DG}}\|z^{\!_N}\|_{L^\infty(\Omega)}\nonumber\\
 &\leq& Ch^{\htf}(N+1)\|u-u^{\!_N}_h\|_{L^2(\Omega)}\|e^{\!_N}_h\|_{\!_{DG}}\nonumber\\
 &\leq& Ch^{\htf}(N+1)\|u-u^{\!_N}_h\|_{L^2(\Omega)}\|e^{\!_N}_h\|_{\!_{DG}}.\nonumber
 \end{eqnarray} }
 
Now, by combining (\ref{bilinear1}) and (\ref{bilinear2}), we have proved the following theorem: 
 \begin{theorem}\label{theorem: dg+}
 Let $u^{\!_N}\in H^{2}(\Omega)$ be the solution of the truncated boundary value problem \emph{(\ref{tdtnbvp})} and $u^{\!_N}_h$ be the solution of the discrete problem~\emph{(\ref{pwdg2})}. Then there exists a constant $C$ depending only on $\Omega$, the flux parameters $\alpha$, $\beta$ and $\delta$ and the shape regularity parameter $\mu$, but independent of $k$,$u^{\!_N}$, $u^{\!_N}_h$, $N$, and $h$ such that 
 \begin{eqnarray}
 \|u^{\!_N}-u^{\!_N}_h\|_{\!_{L^2(\Omega)}} &\leq& Ch^{\hf}(N+1)\inf_{w^{}_h\in PW(\mathscr{T}_h)}\|u^{\!_N} - w^{}_h\|_{\!_{DG^+}}.
 \end{eqnarray}
 \end{theorem} 
{\bf Proof:}\\
By {(\ref{bilinear1})} and {(\ref{bilinear2})} we have 
\begin{eqnarray}
\|u^{\!_N}-u^{\!_N}_h\|^2_{\!_{L^2(\Omega)}} &=& \mathscr{A}_{\!_N}(e^{\!_N}_h,z^{\!_N}-z^c_h)+\mathscr{A}_{\!_N}(e^{\!_N}_h,z^c_h-z^{}_h)\nonumber\\
\end{eqnarray}
The result follows from the error estimate in the $\|\cdot\|_{\!_{DG}}$ norm from Proposition~{\ref{eest1}}.\ep
 
 We need to estimate the term $\displaystyle\inf_{w^{}_h\in PW(\mathscr{T}_h)}\|u^{\!_N} - w^{}_h\|_{\!_{DG^+}}$ that appears in Theorem~\ref{theorem: dg+} along the lines of results due to A. Moiola in \cite{moiolathesis} for the Helmholtz equation with an impedance boundary condition. In \cite{moiolathesis}, detailed error estimates for the approximation of solutions of the homogeneous Helmholtz equation by plane waves are shown. To derive these, two assumptions are made. The first concerns the domain $K\subset\Omega$ (in particular a generalized triangle) where we will approximate a solution of the Helmholtz equation by plane waves (see Lemma 3.1.1 of \cite{moiolathesis}).
The second assumption concerns the distribution of plane wave directions on the unit circle (see Lemma 3.4.3 \cite{moiolathesis}).

 Now assume $w$ is a solution of the homogeneous Helmholtz equation and $w\in H^{m+1}(K)$, where $1\leq m\in\mathbb{Z}$. Assume also $q\geq 2m+1$. Define the $k$-weighted norm $\|w\|_{\!_{\ell,k,K}}$ by
 \begin{eqnarray}
 \|w\|_{\!_{\ell,k,K}}&=& \left(\sum_{j=0}^\ell k^{2(\ell-j)}|w|^2_{\!_{j,K}}\right)^\hf,\;\;\;\forall w\in H^\ell(K),\;\;k>0.
 \end{eqnarray}
 and for every $0\leq j\leq m+1$, consider $\varepsilon_0,\varepsilon_1,\varepsilon_2$ defined in equation (4.20) of \cite{moiolathesis}: 
 \begin{eqnarray}\label{amterm}
 \varepsilon_j:=\left(1+(kh)^{j+6}\right)e^{\left(\frac{7}{4}-\frac{3}{4}\rho  \right)kh}h^{m+1-j}\left[\left(\frac{\log(q+2)}{q}\right)^{({m+1-j})}+\frac{1+(kh)^{q-m+2}}{(c_0(q+1))^{\frac{q}{2}}}\right],\nonumber\\
 \end{eqnarray}
 where $\rho$ is a parameter related to the shape regularity of the elements, such that for a mesh with shape regularity $\mu$, $\rho = (2\mu)^{-1}$. The constant 
 $c_0$ measures the distribution of the directions $\{\bdr_\ell\}$ (Lemma 3.4.3 \cite{moiolathesis}).
  
 Now the general approximation result of Moiola is: 
 \begin{theorem}\emph{(approximation by plane waves: Corollary 3.55 \cite{moiolathesis})}\label{theorem: moiola}\\
 Let $u\in H^{m+1}(K)$ be a solution of the homogeneous Helmholtz equation, where $K\subset\mathbb{R}^2$ is a domain that satisfies assumptions in Lemma 3.1.1 of \cite{moiolathesis}. Suppose the directions $\left\lbrace \bdr_\ell \right\rbrace_{\!_{\ell=-q,\cdots,q}}$ satisfy the assumptions in Lemma 3.4.3 of \cite{moiolathesis}. Then there exists $\vec{\alpha}\in \mathbb{C}^p$ such that for every $0\leq j\leq m+1$,
 \begin{eqnarray}
\|u - \displaystyle\sum_{\ell=1}^p \alpha_{\!_k}e^{ik\bx\cdot\bdr_\ell}\|_{\!_{j,k,K}} &\leq& C\varepsilon_j \|u\|_{\!_{m+1,k,K}}
 \end{eqnarray}
 where $C$ depends on $j,m$ and the shape of $K$.
 \end{theorem}
 
  In the next Lemma, we state best approximation error estimates in the $\|\cdot\|_{\!_{DG^+}}$ norm for the Helmholtz equation with a DtN boundary condition. 
 The error bounds for $\|u-\xi\|_{\!_{0,\cE}}$ and $\|\nabla_h(u-\xi)\|_{\!_{0,\cE}}$ on the skeleton of the mesh can be found in Lemma 4.4.1 of \cite{moiolathesis}. We state them here for completeness.
 For convenience of notation, define $$\gamma:= N^2R^{-2}(1+kR)^2$$

 \begin{lemma}\label{kapita1}
 Assume that the directions $\{\bdr_\ell\}$ satisfy the assumptions of Theorem~\ref{theorem: moiola}.  Given $u\in T(\mathscr{T}_h)\cap H^{m+1}(\Omega)$, $m\geq 1$, $q\geq 2m+1$, there exists $\xi\in PW(\mathscr{T}_h)$ such that we have the following estimates
 \begin{eqnarray}
 \|u-\xi\|^2_{\!_{0,\cE}} &\leq& C\varepsilon_0\left(\varepsilon_0 h^{-1}+\varepsilon_1\right)\|u\|^2_{\!_{m+1,k,\Omega}}\\
 \|\nabla_h(u-\xi)\|^2_{\!_{\cE}} &\leq& C\varepsilon_1\left(\varepsilon_1 h^{-1}+\varepsilon_2\right)\|u\|^2_{\!_{m+1,k,\Omega}}\\
 \left|\Imag \int_{\Gamma_R}\DtN_{\!_N}(u-\xi)\overline{(u-\xi)}\;ds\right| &\leq& Ck\varepsilon_0(\varepsilon_0 h^{-1}+\varepsilon_1)\|u\|^2_{\!_{m+1,k,\Omega}}   \\
 \|\DtN_{\!_N}(u-\xi)\|^2_{\!_{0,\Gamma_{\!_R}}} &\leq& Ck^{-1}\gamma\varepsilon_0\left(\varepsilon_0 h^{-1}+\varepsilon_1 \right)\|u\|^2_{\!_{m+1,k,\Omega}}\\
 \|u-\xi\|^2_{\!_{DG^+}}&\leq& C\left[\varepsilon_0(\varepsilon_0 h^{-1} +\varepsilon_1)(k^{-1}\gamma +k)\right.\nonumber\\
 && + \left.k^{-1}\varepsilon_1(\varepsilon_1 h^{-1}+\varepsilon_2)\right]\|u\|^2_{\!_{m+1,k,\Omega}}
 \end{eqnarray}
 where $C$ is independent of $h,N,k,p,\xi,\{\bdr_\ell\}$ and $u$.
 \end{lemma} 
{\bf Proof:}\\
We have by the trace estimate and Theorem~{{\ref{theorem: moiola}}}
\begin{eqnarray*}
\|u-\xi\|^2_{\!_{0,\partial K}} &\leq& C\left(h_K^{-1}\|u-\xi\|^2_{\!_{0,K}} +\|u-\xi\|_{\!_{0,K}}|u-\xi|_{\!_{1,K}} \right)\\
&\leq& C\varepsilon_0\left(\varepsilon_0+\varepsilon_1\right)\|u\|^2_{\!_{m+1,k,\Omega}},\\
\|\nabla_h(u-\xi)\|^2_{\!_{0,\partial K}} &\leq& C\left(h_K^{-1}|u-\xi|^2_{\!_{1,K}} +|u-\xi|_{\!_{1,K}}|u-\xi|_{\!_{2,K}} \right)\\
&\leq& C\varepsilon_1\left(\varepsilon_1+\varepsilon_2\right)\|u\|^2_{\!_{m+1,k,\Omega}}.
\end{eqnarray*}
Now to prove the error estimate for terms with the DtN map, recall Lemma~\ref{hankel}
\begin{eqnarray*}
\left|\Imag\int_{\Gamma_R} \DtN_{\!_N} v\overline{v}\;ds\right| &=& 4\sum_{|m|\leq N} \frac{|v_m|^2}{|H^{(2)}_m(kR)|^2}\nonumber\\
\mbox{and}&&\\
kR|H^{(2)}_{m}(kR)|^2 &\geq& \frac{2}{\pi},\;\;\;|m|\geq 1.
\end{eqnarray*}
We have,
\begin{eqnarray*}
\left|\Imag \int_{\Gamma_R}\DtN_{\!_N}(u-\xi)\overline{(u-\xi)}\;ds\right| &=& 4\sum_{|m|\le N} \frac{|u_m-\xi_m|^2}{|H^{(2)}_m(kR)|^2}\\
&\leq& 2\sum_{|m|\leq N}\pi kR|u_m-\xi_m|^2\\
&\leq& Ck\|u-\xi\|^2_{\!_{0,\Gamma_{\!_R}}}\\
&\leq& Ck\sum_{K\in\mathscr{T}^R_h}\|u-\xi\|^2_{0,\partial K}.
\end{eqnarray*}
Recalling the definition (\ref{mdef}) of $m_{\!_0}$, and by the bounds in (\ref{hankelpositive}) and (\ref{hankelzero}) 
\begin{eqnarray*}
\|\DtN_{\!_N}(u-\xi)\|^2_{\!_{0,\Gamma_{\!_R}}} &=& 2\pi kR\sum_{|m|\le N}m_{\!_0}^2\left|\frac{1}{m_{0}} \frac{H^{(2)'}_{m}(kR)}{H^{(2)}_{m}(kR)}\right|^2|u_m-\xi_m|^2\\
&\leq& CN^2R^{-2}(1+kR)^2\|u-\xi\|^2_{\!_{0,\Gamma_{\!_R}}}\\
&\leq& CN^2R^{-2}(1+kR)^2\sum_{K\in\mathscr{T}^R_h}\|u-\xi\|^2_{0,\partial K}.
\end{eqnarray*}
To prove the last error estimate, note that
{\begin{eqnarray*}
\|u-\xi\|^2_{\!_{DG^+,N}} &\leq& C\left(k\|u-\xi\|^2_{0,\cE} + k^{-1}\|\nabla_h(u-\xi)\|_{0,\cE}\right. \\ &&\left. +\left|\Imag\int_{\Gamma_R}\DtN_{\!_N}(u-\xi)\overline{(u-\xi)}\;ds\right|\right.\\&&\left. +k^{-1}\|\DtN_{\!_N}(u-\xi)\|^2_{0,\Gamma_R}\right).\nonumber
\end{eqnarray*}}
\ep
 
 The first term in the square brackets of (\ref{amterm}) decays algebraically for large $q$ while the second term decays faster than  exponentially. Therefore for large $q$, we have the following order estimates:
 
 \begin{corollary}\label{coro: kapita2}
  Given $u\in T(\mathscr{T}_h)\cap H^{m+1}(\Omega)$, $m\geq 1$, $q\geq 2m+1$, large enough such that the algebraic term in \emph{(\ref{amterm})} dominates the exponentially decaying term, there exists $\xi\in PW(\mathscr{T}_h)$ such that we have the following estimate  
  \begin{eqnarray*}
  \|u-\xi\|_{\!_{DG^+}} &\leq& CNh^{m-\hf}(h+N^{-1})\;\left(\frac{\log(q+2)}{q}\right)^{m-\hf} \|u\|_{{m+1,k,\Omega}}
  \end{eqnarray*}
   \end{corollary}
  {\bf Proof:} The result follows easily from Lemma~\ref{kapita1} and the inequality $$\|u^{\!_N}\|_{m+1,k,\Omega}\leq C\|u\|_{m+1,k,\Omega}$$ which holds for some $C$ independent of $N\geq N_0$ (see the proof of Theorem 2.5 in \cite{kmdtn}). \ep 
  
 Using Corollary~\ref{coro: kapita2}, we have the following main theorem of this section:
 \begin{theorem}\label{theorem: discr}
 Let $u^{\!_N}$ be the solution of the truncated boundary value problem \emph{(\ref{tdtnbvp})} and $u^{\!_N}_h$ be the computed solution, $m\geq 1$, $q\geq 2m+1$, large enough such that the algebraic term in \emph{(\ref{amterm})} dominates the exponentially decaying term. There exists a constant $C$ that depends on $k$ and $h$ only as an increasing function of their product $kh$, but is independent of $p$, $u^{\!_N}$, $u^{\!_N}_h$ and $N$ such that   
 \begin{eqnarray}\label{errestkap}
 \|u^{\!_N}-u^{\!_N}_h\|_{\!_{L^2(\Omega)}}&\leq& CN h^{m}(h+N^{-1})\left(\frac{\log(q+2)}{q}\right)^{m-\hf}\|u\|_{\!_{m+1,k,\Omega}}.
 \end{eqnarray}
 \end{theorem}
{\bf Proof:}\\
  The result follows from Corollary~{\ref{coro: kapita2}} and Theorem {\ref{eest1}}.
\ep 

{\bf Remark:}  If $k,N,q,R$ are fixed in the error estimate~(\ref{errestkap}), then
$$\|u^{\!_N}-u^{\!_N}_h\|\leq C_2 h^{m}\|u\|_{\!_{m+1,k,\Omega}}$$ where $C_2$ is independent of $h$. \ep

\section{Numerical Implementation}
\label{implm.sec}
Let $N_h=\dim (PW_p(\mathscr{T}_h))$ be the total number of degrees of freedom associated with the PWDG space $PW(\mathscr{T}_h)$. Obviously $N_h=\sum_{K\in\mathscr{T}_h}p_{\!_K}$. The algebraic linear system associated with the DtN-PWDG scheme is 
\begin{eqnarray}
AU&=&F \label{dtn-pwdg-linear}
\end{eqnarray}
where $A\in \mathbb{C}^{N_h\times N_h}$ is the matrix associated with the sesquilinear form  $\mathscr{A}_{\!_N}(\cdot,\cdot)$ ~(\ref{pwdg4}), and $F\in \mathbb{C}^{N_h}$ is the vector associated with the linear functional $\mathscr{L}_h(\cdot)$, and $U\in\mathbb{C}^{N_h}$, the vector of unknown coefficients of the plane waves in $PW(\mathscr{T}_h)$. More precisely,
the discrete solution can be written in terms of plane waves 
\begin{eqnarray}\label{discsol}
{u^N_h} &=& \sum_{K\in\mathscr{T}_h}\sum_{\ell=1}^{p_{\!_K}} u^{\!_K}_\ell \xi_{\ell}^{\!_K} 
\end{eqnarray}
where the coefficients $u^{\!_K}_\ell\in\mathbb{C}$, and the basis functions $\xi^{\!_K}_\ell$ are propagating plane waves 
\[
\xi^{\!_K}_\ell = \begin{cases} \exp(ik\bx\cdot \bdr^{\!_K}_\ell) &\mbox{if } \bx\in K \\ 
0 & \mbox{elsewhere }  \end{cases}
\]
and the directions are given by 
\[
\bdr^{\!_K}_\ell = \left(\cos{\frac{2\pi\ell}{p_{\!_K}}},\sin{\frac{2\pi\ell}{p_{\!_K}}}\right)
\]
Rewriting~(\ref{discsol}) in vector form 
\begin{eqnarray*}
{u^{\!_N}_h} &=& \sum_{j=1}^{N_h} u^h_j \xi^h_j.
\end{eqnarray*}
then 
\begin{eqnarray*}
U &=& \left[u^h_1,\cdots, u^h_{N_h}\right]^T.
\end{eqnarray*}
The global stiffness matrix associated with the sesquilinear form $\mathscr{A}_{\!_N}(u^{\!_N}_h,v_h)$ is  
 $$A = A_{{\small{\!_{Int}}}} + A_{{\small{\!_{Dir}}}} + A_{{\small{\!_{R,loc}}}} + A_{{\small{\!_{DtN}}}}$$ 
 where $A_{{\small{\!_{Int}}}}$ is the contribution from interior edges  $\cEI$, $A_{{\small{\!_{Dir}}}}$ is the contribution from edges on the Dirichlet boundary $\cED$. Terms in $\mathscr{A}_{\!_N}(u^{\!_N}_h,v_h)$ defined on $\cEB$ but with no DtN map contribute $A_{{\small{\!_{R,loc}}}}$. Their computation is standard for PWDG methods.
 
 However the term $A_{{\small{\!_{DtN}}}}$ is computed globally since the DtN map involves an integral on the entire boundary $\Gamma_{\!_R}$, and we give details of the calculation now.
 
To compute the stiffness matrix $A_{{\small{\!_{DtN}}}}$, we introduce the space $W_{\!_N}$ of trigonometric polynomials
\[
W_{\!_N}:=\mbox{span}\left\lbrace e^{i n\theta}: -N\leq n\leq N\right\rbrace.
\]
The projection operator $\Pro_{\!_N}: L^2(\Gamma_{\!_R})\rightarrow W_{\!_N}$ on the artificial boundary $\Gamma_{\!_R}$ onto the $2N+1$ dimensional space of trigonometric polynomials is defined as  
\begin{eqnarray}
\int_{\Gamma_{\!_R}}\left(\Pro_{\!_N} \varphi - \varphi\right)\overline{\eta}\;ds &=& 0,
\end{eqnarray}
where $\varphi\in L^2(\Gamma_{\!_R})$, and $\eta\in W_{\!_N}$. In practice, integrals on $\Gamma_{\!_R}$ are computed elementwise by Gauss-Legendre quadrature using 20 points per edge.

Let $w$ be the projection of the computed solution, $u^{\!_N}_{h}$: $$w := \Pro_{\!_N} u^{\!_N}_h = \sum_{\ell=-N}^N w_{\!_\ell} \eta_{\!_\ell},$$ where $\eta_{\!_\ell}=e^{i\ell\theta}$ is a basis function of $W_{\!_N}$.
Then $$w_{\!_\ell} = \frac{1}{2\pi R}\sum_{j=1}^{N_h} u^h_j\int_{\Gamma_{\!_R}} \xi^h_j \overline{\eta_{\!_\ell}}\;ds = \frac{1}{2\pi R}\sum_{j=1}^{N_h} M_{\!_{\ell j}}u^h_j$$
where the components of the projection matrix $M\in\mathbb{C}^{(2N+1)\times N_h}$ are defined as $$M_{\!_{\ell j}} = \int_{\Gamma_{\!_R}} \xi^h_j\overline{\eta_{\!_\ell}}\;ds.$$
With these coefficients, we can write 
$${W} = \frac{1}{2\pi R} M U,$$
where $U$ is the vector of the unknown coefficients of $u^N_{h}$ previously defined, and 
\begin{eqnarray*}
W &=&\left[w_{\!_{-N}},\cdots, w_{\!_N}\right]^T.
\end{eqnarray*}
Let $v_h = \xi^h_j$. Then $$\Pro_{\!_N} v_h = \Pro_{\!_N} \xi^h_j = \frac{1}{2\pi R}Me_j$$
where $e_j$ is a $N_h\times 1$ vector with one on the $j^{th}$ coordinate, and zero otherwise. Then choosing $v_h=\xi^h_{j}$, the DtN term is computed as 
\begin{eqnarray}
\int_{\Gamma_{\!_R}} \DtN_{\!_N} u^{\!_N}_h \overline{v_h}\;ds &:=& \int_{\Gamma_{\!_R}} (\DtN \Pro_{\!_N}u^{\!_N}_h)(\overline{\Pro_{\!_N}v_h})\;ds\nonumber\\
&=& \frac{1}{2\pi R}\left(M^{\star} T M U\right)_{\!_j},
\end{eqnarray}
$j=1,\cdots,N_h$ and $M^{\star}$ is the conjugate transpose of $M$, and $T$ is the $(2N+1)\times (2N+1)$ diagonal matrix 
$$T = 
\left(\begin{matrix}
\zeta_{\!_{-N}} &  0  & \ldots & 0\\
0  &  \zeta_{\!_{-N+1}} & \ldots & 0\\
\vdots & \vdots & \ddots & \vdots\\
0  &   0       &\ldots & \zeta_{\!_{N}}
\end{matrix}\right)
$$
with diagonal entries $\zeta_{m} = k\frac{H_m^{(2)'}(kR)}{H_m^{(2)}(kR)}$.

To compute terms on $\Gamma_{\!_R}$ involving normal derivatives, observe that 
\begin{eqnarray}
\nabla_h u^{\!_N}_h\cdot \bn &=& \sum_{j=1}^{N_h} u^h_j (ik \bdr_j\cdot \bn)\xi^h_j.
\end{eqnarray}
Let $$z = \Pro_{\!_N}\left(\nabla_h u^{\!_N}_h\cdot\bn\right) = \sum_{m=-N}^N z_{\!_\ell} \eta_{\!_\ell}.$$
We have 
\begin{eqnarray}
2\pi Rz_{\!_\ell} &=& \int_{\Gamma_{\!_R}}\Pro_{\!_N} \left(\sum_{j=1}^{N_h} u^h_j(ik\bdr_j\cdot \bn)\xi^h_j\right)\overline{\eta_{\!_\ell}}\;ds\nonumber\\
&=& \sum_{j=1}^{N_h} u^h_j \int_{\Gamma_{\!_R}}(ik\bdr_j\cdot\bn)\xi^h_j \overline{\eta_{\!_\ell}}\;ds.
\end{eqnarray}
Then 
$$Z = \frac{1}{2\pi R}M_{\!_D}U,$$
where 
\begin{eqnarray*}
Z &=& \left[z_{\!_{-N}},\cdots, z_{\!_{N}}\right]^T
\end{eqnarray*}
and $M_{\!_D}$ is the $(2N+1)\times N_h$ projection-differentiation matrix with entries
$$({M_{\!_D}})_{\!_{\ell j}} = \int_{\Gamma_{\!_R}}(ik\bdr_j\cdot\bn)\xi^h_j \overline{\eta_{\!_\ell}}\;ds.$$

Choosing $v_h = \xi^h_{\!_j}$, we have
\begin{eqnarray}\label{dtn2}
&& -\frac{1}{ik}\int_{\Gamma_{\!_R}} \delta(\nabla_hu^{\!_N}_h\cdot\bn-\DtN_{\!_N} u^{\!_N}_h)\overline{(\nabla_h \xi^h_j\cdot \bn-\DtN_{\!_N} \xi^h_j)}\;ds\nonumber\\
&=& -\frac{1}{ik}\int_{\Gamma_{\!_R}}\delta(\nabla_h u^{\!_N}_h\cdot\bn)\overline{( \nabla_h \xi^h_j\cdot\bn)}\;ds +\frac{1}{ik}\int_{\Gamma_{\!_R}} \delta (\DtN_{\!_N} u^{\!_N}_h)\overline{(\nabla_h \xi^h_j\cdot\bn)}\;ds\nonumber\\
&& +\frac{1}{ik}\int_{\Gamma_{\!_R}}\delta(\nabla_h u^{\!_N}_h\cdot\bn)\overline{(\DtN_{\!_N} \xi^h_j)}\;ds
 -\frac{1}{ik}\int_{\Gamma_{\!_R}}\delta (\DtN_{\!_N} u^{\!_N}_h)\overline{(\DtN_{\!_N} \xi^h_j)}\;ds\nonumber\\
 &=& I_1 + I_2 + I_3 +I_4.
\end{eqnarray}
The first term $I_1$ is computed locally on each edge, in the standard way. The second term $I_2$ is computed as follows
\begin{eqnarray}\label{term1}
I_2 &=& \frac{1}{ik}\int_{\Gamma_{\!_R}}\delta \left(\DtN \Pro_{\!_N} u^{\!_N}_h\right)\overline{\left(\nabla_h \xi^h_{\!_j}\cdot\bn\right)}\;ds\nonumber\\
&=& \frac{\delta}{2\pi ik R}[M_{\!_D}^{\star} T M U]_j.
\end{eqnarray}
The third term is computed as 
\begin{eqnarray}\label{term2}
I_3 &=& \frac{1}{ik}\int_{\Gamma_{\!_R}}\delta\left(\nabla_h u^{\!_N}_h\cdot\bn\right)\left(\overline{\DtN \Pro_{\!_N}\xi^h_j}\right)\;ds\nonumber\\
&=& \frac{\delta}{2\pi ik R}[(TM)^{\star}M_{\!_D}U]_j.
\end{eqnarray}
To compute the fourth term,
\begin{eqnarray}
I_4 &=& -\frac{1}{ik}\int_{\Gamma_{\!_R}}\delta (\DtN \Pro_{\!_N}u^{\!_N}_h)(\overline{\DtN \Pro_{\!_N} \xi^h_j})\;ds\nonumber\\&=& -\frac{\delta}{2\pi ik R}[(TM)^{\star}(TM)U]_j.
\end{eqnarray}
It follows that the contribution to the global stiffness matrix from terms involving the DtN map on $\Gamma_{\!_R}$ is the $N_h\times N_h$ matrix
\begin{eqnarray}\label{dtnmat}
A_{{\small{DtN}}} &=& -\frac{1}{2\pi R}M^{\star}TM +\frac{\delta}{2 ik\pi R}\left[M_{\!_D}^{\star}TM + (TM)^{\star}M_{\!_D} - (TM)^{\star}(TM)\right].\nonumber\\
\end{eqnarray}
\section{Numerical Results}
\label{results.sec}
In this section, we numerically investigate convergence of the DtN-PWDG scheme. We consider the scattering of acoustic waves by a sound-soft obstacle, as modeled by the boundary value problem (\ref{helmholtzbc}). In our numerical experiments, we consider both the impedance boundary condition (IP-PWDG)
for the scattered field
\begin{eqnarray}
&&\frac{\partial u}{\partial \bn} +iku = 0,\;\;\mbox{on}\;\Gamma_{\!_R}
\end{eqnarray}
and the DtN boundary condition,
\begin{eqnarray}
&&\frac{\partial u}{\partial \bn} - \DtN_{\!_N}u = 0,\;\;\mbox{on}\;\Gamma_{\!_R}.
\end{eqnarray}
In all numerical experiments, the Dirichlet boundary condition is imposed on the scatterer
\begin{eqnarray}
&&u = -u^{\small{\mbox{inc}}},\;\;\mbox{on}\;\;\Gamma_{\!_D}.
\end{eqnarray}
where $u^{\small{\mbox{inc}}}= e^{ik\bx\cdot\bdr}$ is a plane wave incident field propagating in the direction $\bdr$ relative to a negative orientation of the $xy$ axis, due to our choice of the time convention in the definition of the time-harmonic field. The computational domain is the annular region between $\Gamma_{\!_R}$ and $\Gamma_{\!_D}$.

For scattering from a disk, the scatterer $D$ is a circle of radius $a = 0.5$, while the artificial boundary on which the DtN map is imposed is a circle of radius $R = 1$, centered at the origin. All computations are done in MATLAB. The code used in our numerical experiments is based on the 2D Finite Element toolbox LEHRFEM \cite{lehrfem}. The outer edges on $r=a$, $r=R$ are parametrized in polar coordinates, and high order Gauss-Legendre quadrature (20 points per edge) is used to compute all integrals defined on edges on $\Gamma_{\!_R}$ and $\Gamma_{\!_D}$. Curved edges are used in order to eliminate errors that arise from using an approximate polygonal domain.

\begin{figure}
\begin{center}
\resizebox{0.40\textwidth}{!}{\includegraphics{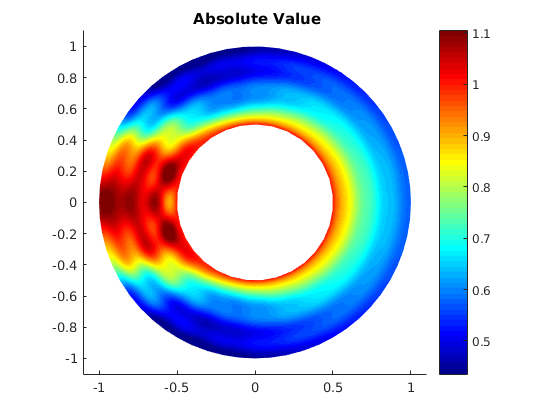}}
\resizebox{0.40\textwidth}{!}{\includegraphics{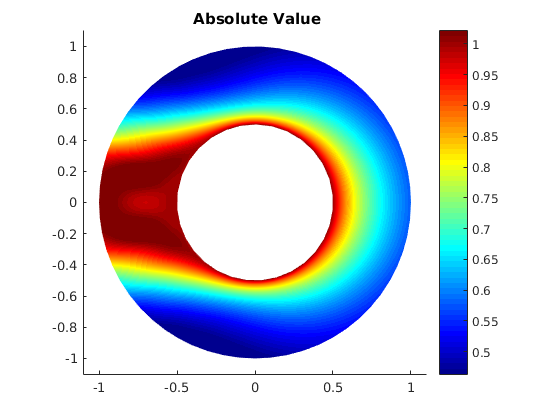}}
\resizebox{0.40\textwidth}{!}{\includegraphics{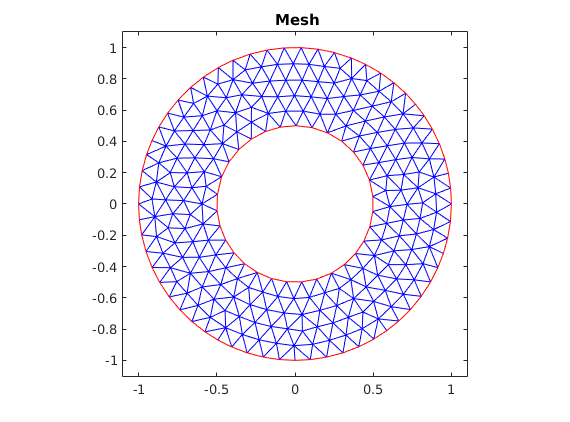}}
\resizebox{0.40\textwidth}{!}{\includegraphics{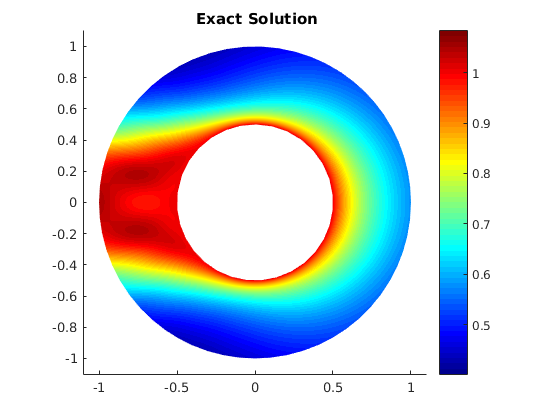}}
\resizebox{0.40\textwidth}{!}{\includegraphics{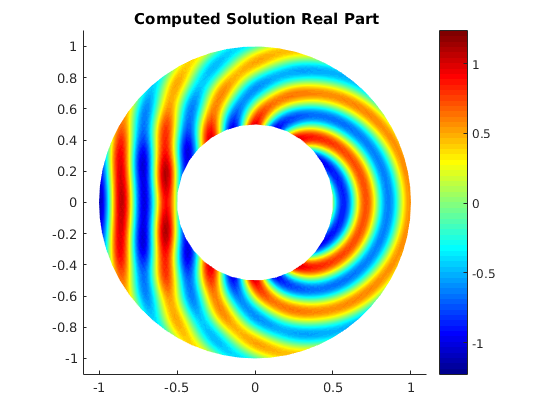}}
\resizebox{0.40\textwidth}{!}{\includegraphics{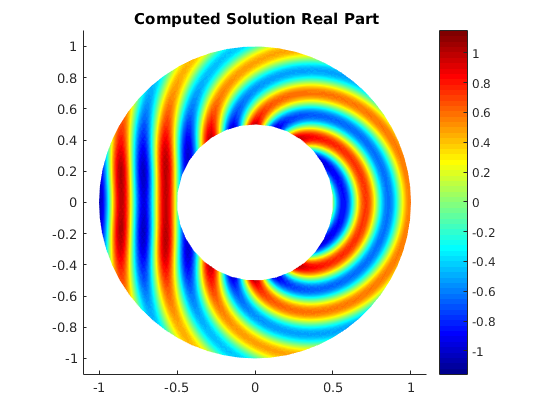}}
\caption{Scattering from a sound-soft disk: $a=0.5$, $R=1$, $p=15$ plane waves per element, $k=7\pi$, $N=20$ Hankel functions. Top left: absolute value of the solution computed using impedance boundary conditions. Top right: absolute value of the solution using DtN boundary conditions. Middle left: the mesh. Middle right: absolute value of the exact solution. Bottom left: real part computed using the impedance boundary conditions. Bottom right: real part computed using the DtN boundary conditions. }\label{fig:scatcirc1}
\end{center}
\end{figure}

 The exact solution of the scattering problem (\ref{helmholtzbc}), in the case of a circular scatterer, is given in polar coordinates by (see Section 6.4 of Colton \cite{coltonpde})
 \begin{eqnarray}\label{eqn4: exactscat}
 u(r,\theta) &=& -\left[\frac{J_0(ka)}{H^{(2)}_0(ka)}H^{(2)}_0(kr)+2\sum_{m=1}^\infty i^m\frac{J_m(ka)}{H^{(2)}_m(ka)}H^{(2)}_m(kr)\cos{m\theta}\right]\nonumber \\
\end{eqnarray}
For numerical experiments, we take $N=100$ to truncate the exact solution. This value of $N$ is sufficient for the wavenumbers considered.
\paragraph*{Experiment 1: Scattering from a disk}
Our main example is a detailed investigation of the problem of scattering of a plane wave from a disk. We choose the DtN or impedance boundary $\Gamma_{\!_R}$ to be concentric (this improves the accuracy of the impedance boundary condition). In Fig~\ref{fig:scatcirc1}, we compare density plots of the approximate solution by both methods using the same discrete PWDG space.
Comparison of the shadow region of the impedance boundary condition solution in Fig~\ref{fig:scatcirc1} (top left) with that of the exact solution demonstrates the effect of spurious reflections from the artificial boundary. The DtN-PWDG solution in the top right panel of Fig~\ref{fig:scatcirc1} shows greater fidelity with the exact solution. There is little difference between the shadow regions of the DtN solution and the exact solution. It is interesting to note that from the plots of the real parts of the solution in the bottom row of Fig~\ref{fig:scatcirc1}, there is little obvious difference between the solutions. However our upcoming and more detailed analysis shows significant improvements from the DtN boundary condition.

\begin{figure}
\begin{center}
\resizebox{0.45\textwidth}{!}{\includegraphics{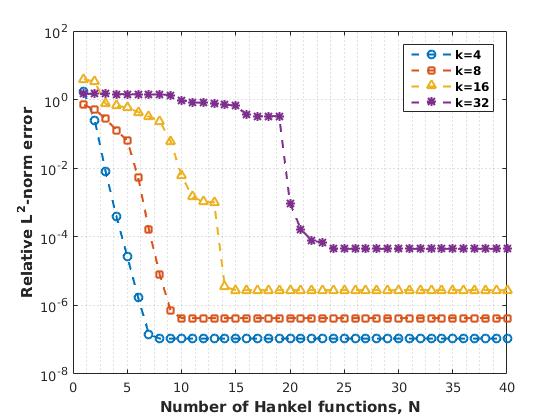}}
\resizebox{0.45\textwidth}{!}{\includegraphics{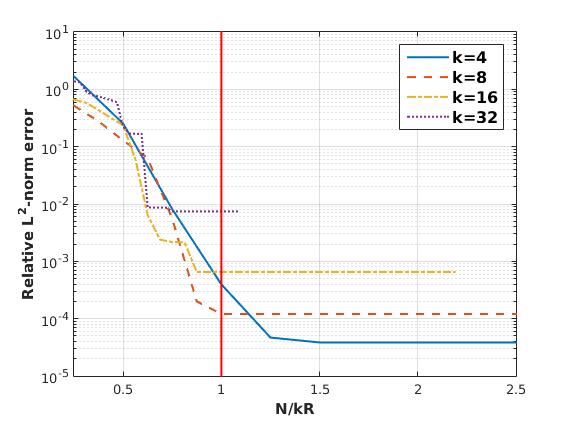}}
\resizebox{0.45\textwidth}{!}{\includegraphics{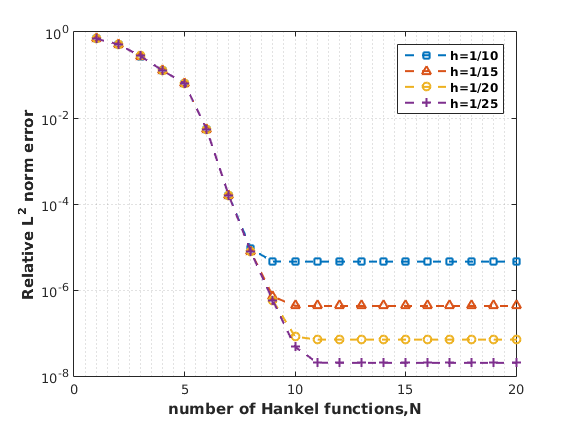}}
\caption{Scattering from a disk: Top: semi-log plot of the relative $L^2$-norm error vs maximum order $N$ number of the Hankel functions in the DtN expansion, for $k=4,8,16,32$, $p=11$, $h=1/15$. Middle: log of the relative $L^2$-norm error vs $N/kR$, $p=7$, $h=0.1$. Bottom: log of the relative $L^2$-norm error vs $N$, for $k=8$, $p=11$, $h=1/15$. 
}\label{fig:nversion}
\end{center}
\end{figure}
Our first detailed study investigates the error due to truncation for the DtN-PWDG. We fix a grid and PWDG space and vary $N$ for several wavenumbers $k$. Results are shown in Fig~\ref{fig:nversion}.

Results from the top panel of Fig~\ref{fig:nversion} suggest that for all values of $k$ considered there exists an $N_{0,k}$ such that no further improvement in accuracy is possible for $N>N_{0,k}$, for a fixed number of plane waves $p$ and a fixed mesh width $h$. Taking $N>N_{0,k}$ does not improve the accuracy of the solution. The error is then due to the PWDG solution. There are three phases in the plots: (i) a pre-convergence phase, where increasing $N$ has little effect on convergence (ii) convergence phase where rapid exponential convergence of the relative error with respect to $N$ is observed. (iii) post-convergence phase when $N>N_{0,k}$ and optimal convergence has been reached. The middle plot suggests that $N_{0,k}\geq 1.2\;kR$ is sufficient to reach the optimal order of accuracy, which agrees with the rough numerical rule of thumb $N_{0,k}> kR$. In the bottom plot, we note that the exponential rate of convergence in the convergence phase is independent of $h$, however the final $N_{0,k}$ and  error depend on $h$.

\begin{figure}
\begin{center}
\resizebox{0.40\textwidth}{!}{\includegraphics{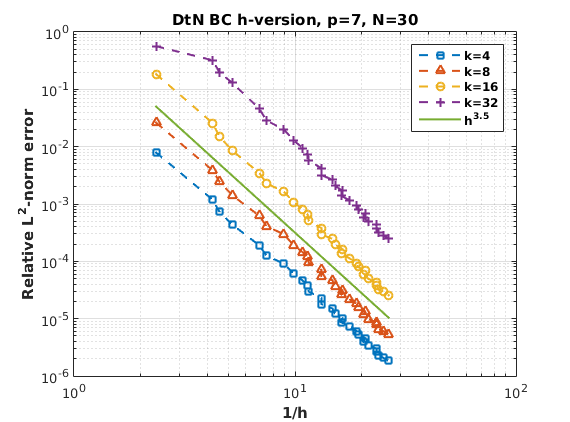}}
\resizebox{0.40\textwidth}{!}{\includegraphics{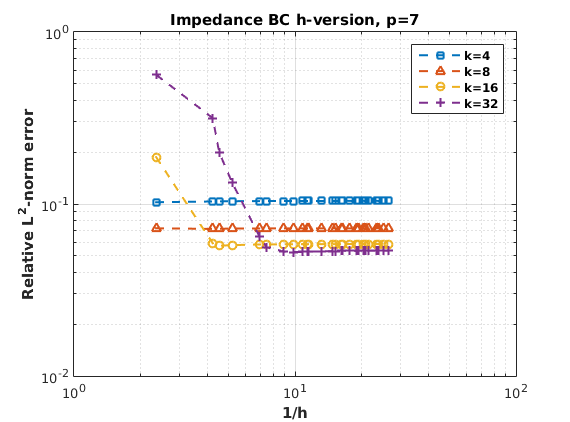}}
\caption{Scattering from a disk: log-log plot of the relative $L^2$-norm error vs $1/h$. Top: DtN-PWDG with $N=30$, $p=7$ plane waves per element. Bottom: IP-PWDG, $p=7$.}\label{fig:hversion}
\end{center}
\end{figure}
Next, we fix $N=30$ (sufficiently large on the basis of the previous numerical results that the error due to the truncation of the DtN map is negligible) and examine $h$-convergence for the DtN-PWDG and standard PWDG with an impedance boundary condition. Results are shown in Fig~\ref{fig:hversion}.

 The top graph suggests that the rate convergence of the relative error with respect to $h$ is independent of $k$ since all curves are roughly parallel. For all $k$ considered, the rate of convergence for $p=7$ is roughly the rate of $3.5$ which exceeds the rate of about $3.0$ predicted in Theorem \ref{theorem: discr}. The actual relative error however, depends on $k$, as expected from consideration of dispersion: higher values of $k$ result in more dispersion error which is unavoidable even in the PWDG method~\cite{GH14}. The bottom plot is for PWDG with an impedance boundary condition. It suggests limited convergence of the relative error computed using impedance boundary conditions. However, the accuracy of the solution increases with the wavenumber, in a way contrasting with the results of the top graph for DtN-PWDG. This suggests that the errors due to the approximate boundary condition exceed those due to numerical dispersion in this example.

\begin{figure}
\begin{center}
\resizebox{0.45\textwidth}{!}{\includegraphics{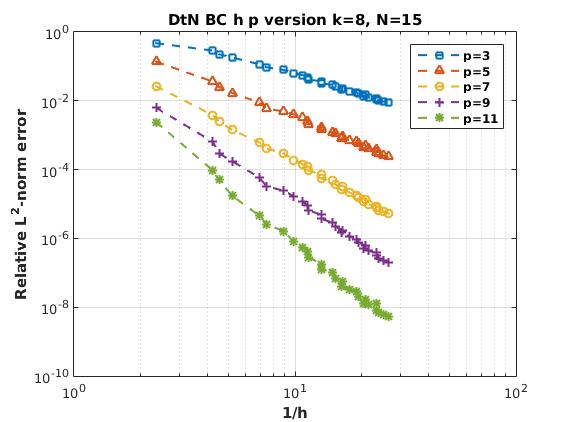}}
\resizebox{0.45\textwidth}{!}{\includegraphics{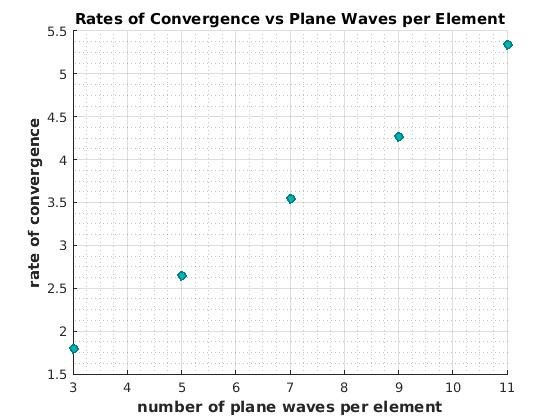}}

\caption{Scattering from a disk:  Top: log-log plot of the relative $L^2$ error vs $1/h$. Bottom: empirical rates of $h$-convergence for different values of $p$. }\label{fig:hpversion}
\end{center}
\end{figure}
Our next example examines $h$ convergence for different choices of $p$. We only consider the DtN-PWDG because of the adverse error characteristics the impedance PWDG shown in Fig~\ref{fig:hversion}. Results for the $h$ and $p$ study are shown in Fig~\ref{fig:hpversion}. 

The results in Fig~\ref{fig:hpversion} top panel show the increased rate of convergence of the PWDG when $p$ is increased. This is clarified in the lower panel. As expected from Theorem~\ref{theorem: discr}, increasing the number of directions of the plane waves per element results in a progressively higher order scheme.
\begin{figure}
\begin{center}
\resizebox{0.45\textwidth}{!}{\includegraphics{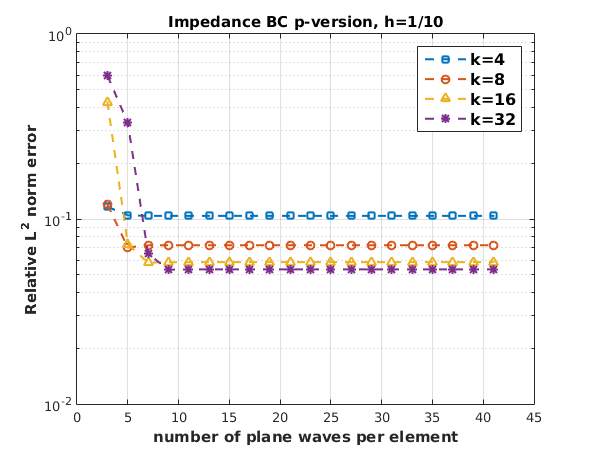}}
\resizebox{0.45\textwidth}{!}{\includegraphics{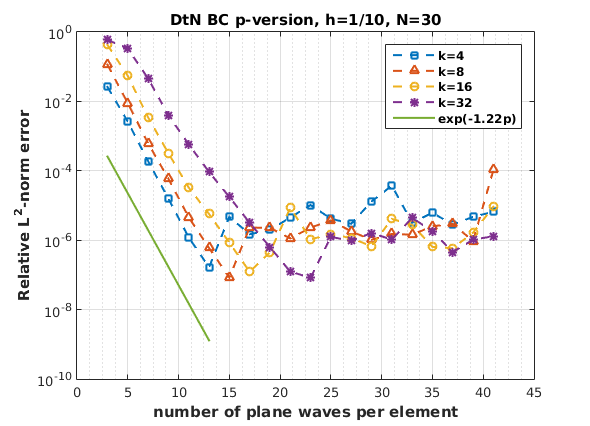}}
\caption{Scattering from a disk: log of the relative $L^2$ error vs $p$ the number of plane waves per element. Top: impedance boundary condition. Bottom: DtN boundary condition with $N=30$, $h=0.1$. }\label{fig:pversion}
\end{center}
\end{figure}

Our final numerical study for scattering from a disk examines $p$ convergence of the DtN-PWDG and impedance PWDG. Results are shown in Fig~\ref{fig:pversion} where we fix $N$ and the mesh size $h$. 

From Fig~\ref{fig:pversion}, as in the case of $h$-convergence, the impedance boundary condition shows limited convergence up to a relative error of about $10\%$, again suggesting that the error due to the boundary condition dominates the error due to the PWDG method in this example. 
The relative $L^2$ error for the DtN boundary condition converges exponentially fast with respect to $p$, the number of plane wave directions per element. However convergence stops due to numerical instability caused by ill-conditioning at a relative error of $10^{-4}\%$. From these experiments we note that the critical number of plane waves needed before numerical instability sets in depends on the wavenumber.

\begin{figure}
\begin{center}
\resizebox{0.40\textwidth}{!}{\includegraphics{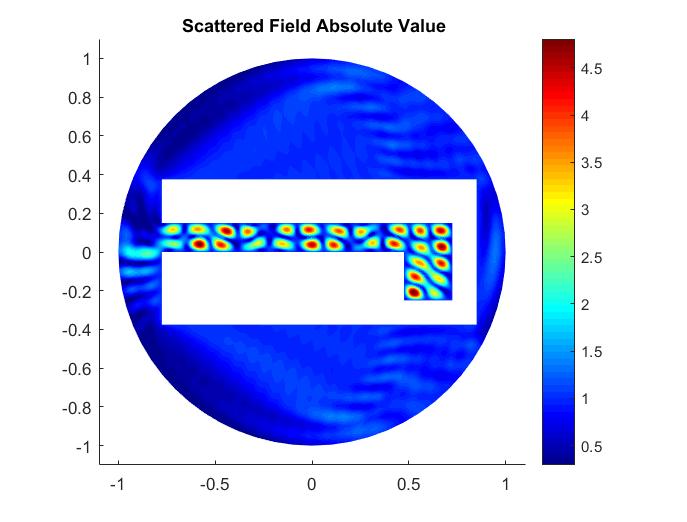}}
\resizebox{0.4\textwidth}{!}{\includegraphics{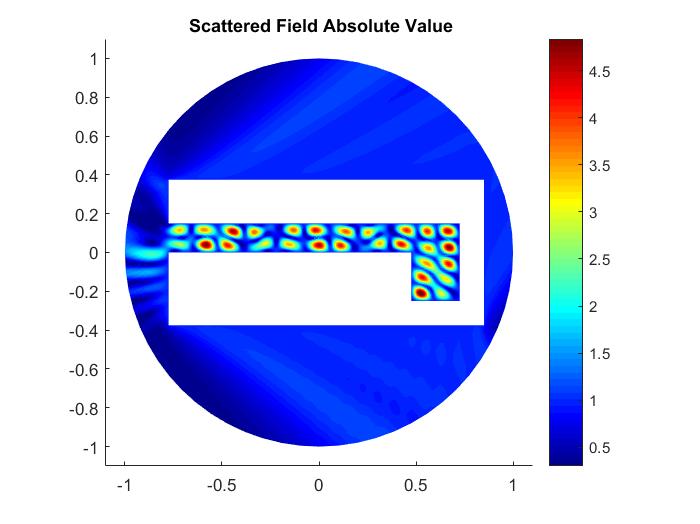}}
\resizebox{0.4\textwidth}{!}{\includegraphics{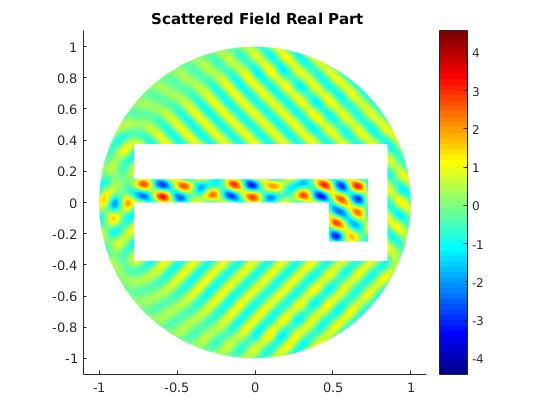}}
\resizebox{0.4\textwidth}{!}{\includegraphics{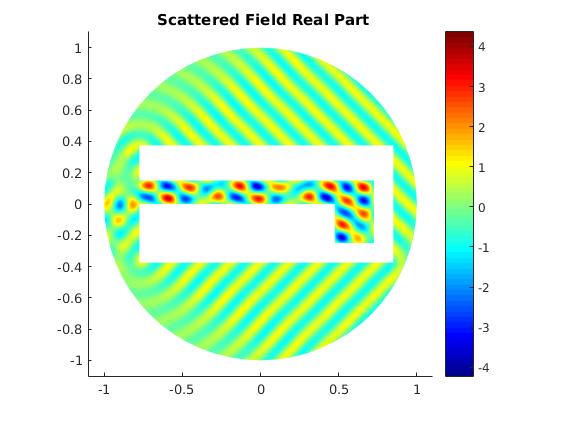}}
\resizebox{0.4\textwidth}{!}{\includegraphics{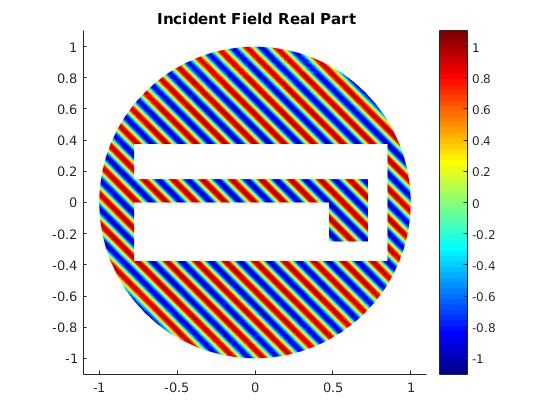}}
\resizebox{0.4\textwidth}{!}{\includegraphics{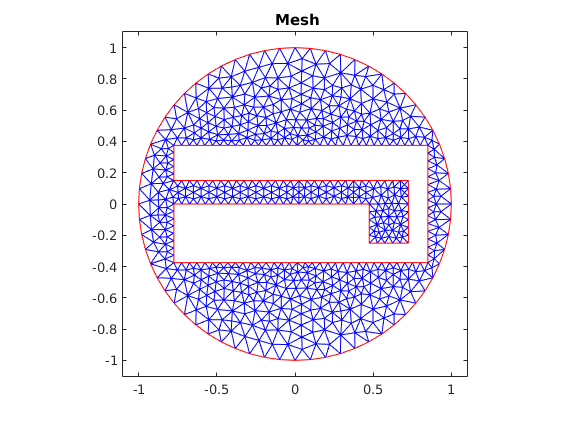}}

\caption{Scattering from a domain with an $L$-shaped cavity, $p=15$ plane waves per element, $k=15\pi$. Top left: absolute value of the scattered field, IP-PWDG. Top right: absolute value of the scattered field, DtN-PWDG. Middle left: real part  of the scattered field, IP-PWDG. Middle right: real part of the scattered field, DtN-PWDG.  Bottom left: The incident field in the direction $\bdr = -\frac{1}{\sqrt{2}}(1\;\;1)$. Bottom right: mesh}
\label{fig:lshape}
\end{center}
\end{figure}
\paragraph*{Experiment 2: Scattering from a resonant cavity}
In our next experiment in Fig~\ref{fig:lshape} we show results for a resonant $L$-shaped cavity. This domain does not satisfy the geometric constraint that the scatterer is star-shaped with respect to the origin. The domain can be included in our theory except we can no longer state $k$-dependent continuity and error estimates. The solution will still be in $H^{\htf+s}(\Omega)$ for some $s>0$. We consider scattering of a plane wave $e^{ik\bx\cdot\bdr}$ from a non-convex domain with an $L$-shaped cavity in the interior, where the direction of propagation of the plane wave is $\bdr = -\frac{1}{\sqrt{2}}(1\;\;1)$. The top left and right panels of Fig~\ref{fig:lshape} show the absolute value of the scattered field computed using IP-PWDG and DtN-PWDG respectively. The IP-PWDG results show reflections on the right hand side of the domain $\Omega$ reminiscent of the poor results in the shadow region for scattering from the disk. These reflections are not visible in the shadow region of the DtN-PWDG solution. 

\begin{figure}
\begin{center}
\resizebox{0.4\textwidth}{!}{\includegraphics{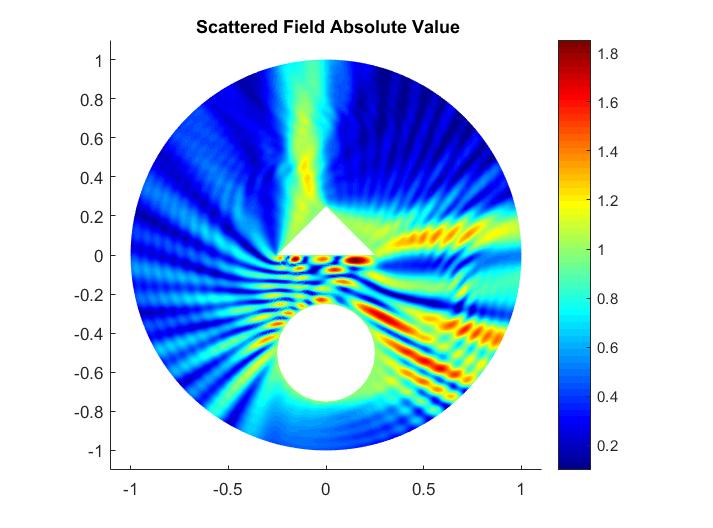}}
\resizebox{0.4\textwidth}{!}{\includegraphics{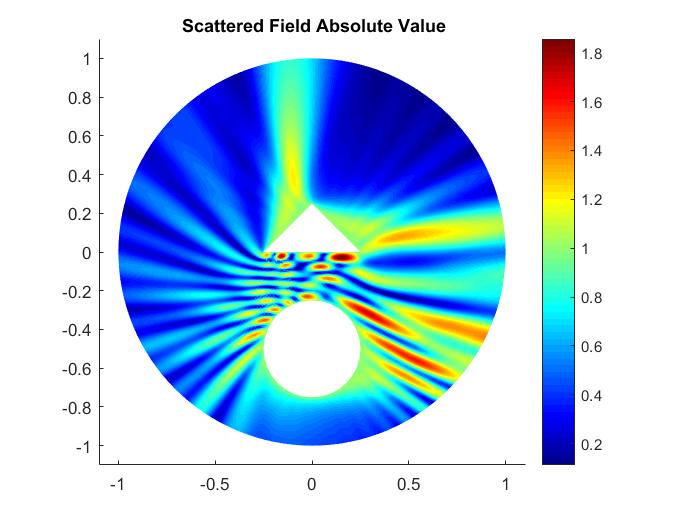}}
\resizebox{0.4\textwidth}{!}{\includegraphics{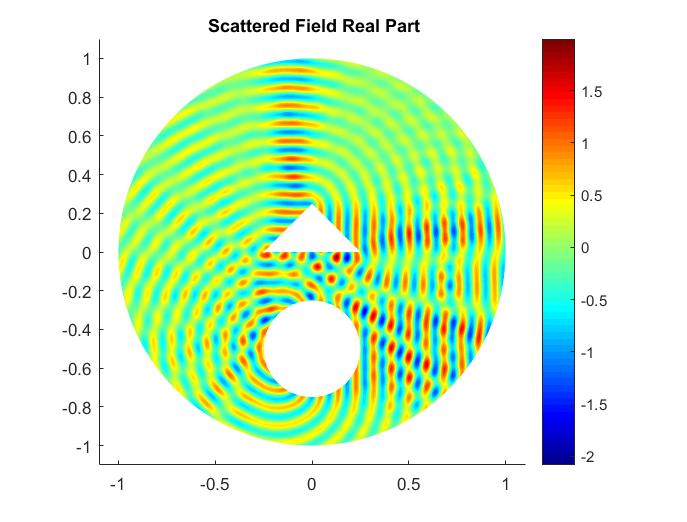}}
\resizebox{0.4\textwidth}{!}{\includegraphics{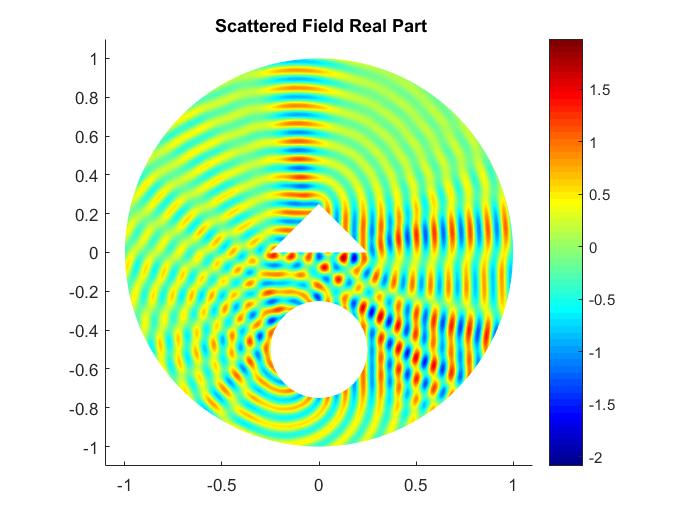}}
\resizebox{0.4\textwidth}{!}{\includegraphics{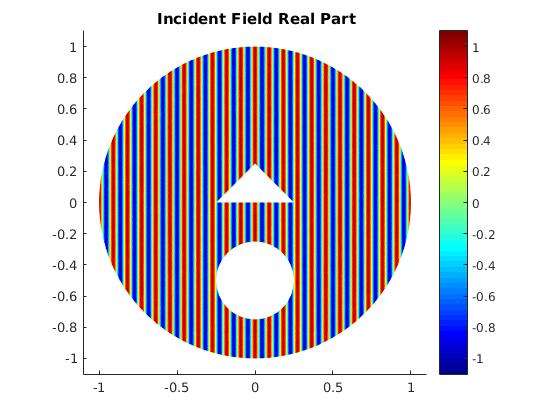}}
\resizebox{0.4\textwidth}{!}{\includegraphics{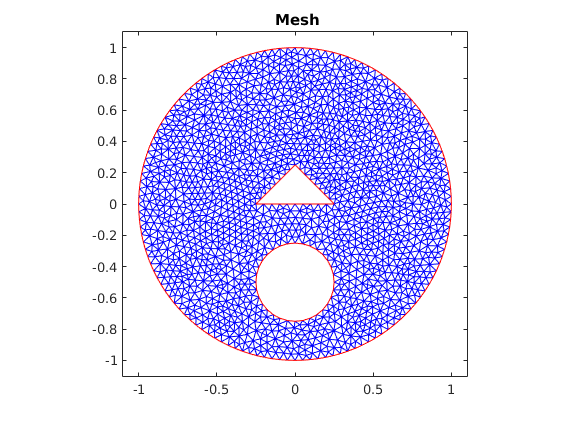}}
\caption{Scattering from a disconnected domain: $k=22\pi$, $p=15$ plane waves per element. Top left: absolute value of the scattered field, IP-PWDG. Top right: absolute value of the scattered field, DtN-PWDG. Middle left: real part of the scattered field, IP-PWDG. Middle right: real part of the scattered field, DtN-PWDG. Bottom left: incident field in the direction $\bdr = (-1\;\;0)$. Bottom right: mesh.}\label{fig:compldtn}
\end{center}
\end{figure}
\paragraph*{Experiment 3: Scattering from a disconnected domain}
In the final set of experiments, we consider scattering of a plane wave incident field from a disconnected domain. The top left solution in Fig~\ref{fig:compldtn} computed using IP-PWDG shows the effect of spurious reflections compared with the smoother shadow region in the DtN-PWDG solution reminiscent of the results in Experiment 1 for scattering from a disk. The real parts of the scattered field are identical to the eye.

This experiment demonstrates the importance of an improved treatment of the absorbing boundary condition for disconnected scatterers.

\section{Conclusion}
\label{conclusion.sec}
We have shown how to couple the PWDG scheme to a non-local boundary condition in order to provide accurate truncation for the scattering problem. Both convergence analysis and numerical results demonstrate good accuracy of the method. Future work would be to use boundary integral equations to approximate the DtN equation. Our formulation should generalize to that case.
\section{Acknowledgements}
The work of the authors was supported in part by NSF grant number DMS-1216620 and AFOSR Grant number FA-9550-11-1-0199. The authors acknowledge the support of the IMA, University of Minnesota during the special year “Mathematics and Optics”.
\bibliographystyle{abbrv}
\bibliography{dtn}

\end{document}